\documentclass[11pt]{article}
 \usepackage{amsmath,amsfonts,amssymb,mathrsfs,cases,mathdots,colordvi,url}
 \usepackage{datetime}
 \usepackage{color}
 \allowdisplaybreaks

 \newtheorem{theorem}{Theorem}[section]
 \newtheorem{proposition}{Proposition}[section]
 
 \newtheorem{remark}{Remark}[section]
 \newtheorem{definition}{Definition}[section]
 \newtheorem{lemma}{Lemma}[section]
 \newtheorem{corollary}{Corollary}[section]
 \newtheorem{example}{Example}[section]

 \renewcommand{\Re}{\mathbb{R}}

 \def \dist{{\rm dist} }

 \def\beginproof{\par\noindent {\bf Proof.}\ \ }
 \def\endproof{\hskip .5cm $\Box$ \vskip .5cm}

 \newcommand{\K}{{\cal K}}

\date{}

\begin{document}

 \title{Exact formula for the second-order tangent set of
the second-order cone complementarity set}

 \author{Jein-Shan Chen\thanks{Department of Mathematics,
 National Taiwan Normal University,
 Taipei 11677, Taiwan, e-mail: jschen@math.ntnu.edu.tw. The author's work is
 supported by Ministry of Science and Technology, Taiwan.} \ \ \ \ \ Jane J. Ye\thanks{Corresponding author. Department of Mathematics
and Statistics, University of Victoria, Victoria, B.C., Canada V8W 2Y2, e-mail: janeye@uvic.ca. The research of this author is
supported by NSERC.}\ \ \ \ \ Jin Zhang\thanks{Department of Mathematics, Southern University of Science and Technology, Shenzhen, P.R. China.
e-mail: zhangj9@sustech.edu.cn.  This author's work is supported by NSFC (11601458, 11871269).} \ \ \ \ \ Jinchuan Zhou\thanks{Department of Statistics, School of Mathematics and Statistics, Shandong University of Technology,
 Zibo 255049, P.R. China, e-mail: jinchuanzhou@163.com. This author's work is supported by National Natural Science Foundation of China (11771255, 11801325, 71871009) and Shandong Province Natural Science Foundation (ZR2016AM07).}}

 \maketitle

 \noindent
 {\bf Abstract.}\
 The second-order tangent set is an important concept in describing the curvature of the set involved. Due to the existence of the complementarity condition, the second-order cone (SOC) complementarity set is a nonconvex set. Moreover, unlike the vector complementarity set, the SOC complementarity set  is not even the union of finitely many polyhedral convex sets. Despite these difficulties,  we succeed in  showing that like the vector complementarity set, the SOC complementarity set is  second-order directionally differentiable and an exact  formula for the second-order tangent set of the SOC complementarity set can be given.  We derive these results  by establishing   the relationship between the second-order tangent set of the SOC complementarity set and the second-order directional derivative of the projection operator over {the} second-order cone, and calculating the  second-order directional derivative of the projection operator over {the} second-order cone. {As an application, we derive second-order necessary optimality conditions for the mathematical program with second-order cone complementarity constraints.}



 \section{Introduction}
In optimization, an important issue is how to approximate the feasible region using derivatives of the function and  the tangent cone of the set involved. Such needs arise in optimality conditions, constraint qualifications and stability analysis when the problem data are perturbed.
  In the same way that second-order derivatives provide quadratic approximations  whereas first-order derivatives only provide linear approximation to a given function, second-order tangent sets provide  better approximation than tangent cones to a set at a point, in particular when the given set is not a polyhedral set or the union of finitely many polyhedral sets. As a result, the second-order tangent sets have been used successfully in second-order optimality conditions, stability analysis,
   and metric subregularity (see e.g. \cite{BCS99,BR05,BS,Cons06,Con16,GM11,GJN10,JN04} {and references therein}). More recently, Gfrerer and Mordukhovich \cite{GM16} use the second-order tangent set to give an estimate of the upper curvature of a set, which is used to study the Robinson regularity of parametric constraint systems.

 In optimization, one often has to deal with a feasible region in the form
  $C:=\{x|\, F(x)\in \Theta\},$
  where $F:\Re^n\rightarrow \Re^m$ is a second-order continuously differentiable  mapping and $\Theta$ is a closed set in $\Re^m$.
 By \cite[Proposition 13.13]{RW98}, under a constraint qualification, the second-order tangent set of the feasible region
 $C$ can be characterized as
 \begin{eqnarray}
 \left .
 \begin{array}{r}
 d\in T_C( x) \\
 w\in T^2_C( x,d)
 \end{array}
 \right \} \quad \Longleftrightarrow \quad \left \{
 \begin{array}{l}
 {\nabla} F(x) d\in T_{\Theta}(F( x)) \\
 {\nabla} F( x) w+d^T {\nabla}^2 F( x)d\in T^2_{\Theta}(F( x);\nabla F( x)d),
 \end{array}
 \right. \label{chainruletangentset}
 \end{eqnarray}
 where $T_C, T_C^2$ denote the tangent cone and the second-order tangent set, respectively (see Definition \ref{sotangent}).
  In the case when $\Theta=\Re_-^{m_1}\times \{0\}^{m_2}$, $m_1+m_2=m$, the system is described by  inequality and equality constraints.
In this  case, since the set $\Theta$ is  polyhedral,     the second-order tangent set of $\Re_-^{m_1}\times \{0\}^{m_2}$ is a polyhedral set, and hence the second-order tangent set of the feasible region is a system of equalities and inequalities involving the second-order derivatives of the constraint mapping $F$ (see, e.g., Bonnans and Shapiro \cite[Formula (3.81)]{BS}), provided a constraint qualification holds.
  In recent years, the second-order cone programming (SOCP) has attracted much attention due to a broad range of applications in fields from engineering, control and finance to robust optimization and combinatorial optimization (see  e.g., \cite{AG} for introduction to the theories and its applications).

Consider the second-order cone defined as
\[
 \K :=  \{ (x_1,x_2) \in \Re \times \Re^{m-1} \, | \,  \| x_2 \| \leq x_1 \},
\]
 where $\| \cdot \|$ denotes the Euclidean norm. Bonnans and Ram\'{i}rez gave the characterization for the second-order tangent set \cite[Lemma 27]{BR05}, and using it to formulate  second-order necessary and sufficient optimality conditions for nonlinear SOCPs.  Since the second-order cone is not polyhedral, the second-order tangent set is not polyhedral \cite{BR05}.

 In recent years,  there are more and more researches on the second-order cone (SOC) complementarity system defined as
\begin{eqnarray*}
&&{\cal K}  \ni G(z)  \perp H(z)  \in {\cal K},
\end{eqnarray*}
 where $u\perp v$ means the vectors $u$ and $v$ are perpendicular, $G(z), H(z): \Re^n\rightarrow \Re^m$. One of the sources of the SOC complementarity system is the Karush-Kuhn-Tucker (KKT) optimality condition  for the second-order cone programming (see e.g. \cite{AG,CT}), and the other is the equilibrium system for a Nash game where the constraints involving second-order cones (see e.g. \cite{HYF}).
We call the closed cone
$$ \Omega
 := \{(x,y) \in \Re^{2m} \, | \, \K \ni x \perp y \in \K\},$$
  the SOC complementarity set (or the complementarity set associated with the second-order cone,
c.f.  \cite{JLZ15}). Using the SOC complementarity set, the SOC complementarity system can be reformulated as
$ (G(z),H(z)) \in \Omega.$ Due to the existence of the complementarity condition, the SOC complementarity set is a nonconvex set. Moreover,   due to the nonpolyhedral structure of the second-order cone ${\cal K}$, the SOC complementarity set is also nonpolyhedral.
 Hence the SOC complementarity set is a difficult object to study in the variational analysis.

 The main goal of this paper is to provide a precise formula for the second-order tangent set to the SOC complementarity set $\Omega$. The projection operator over the second-order cone $\Pi_{\cal K}(x):=\arg\min_{ x'\in {\cal K}}\|x'-x\|$ is one of our main tools in the subsequent analysis.
It is well-known that  the metric projection operator $\Pi_{\cal K}(x)$ provides an alternative characterization of the SOC complementarity set:
 \begin{eqnarray} \label{Omeganew}
 (x,y)\in \Omega
\Longleftrightarrow   \Pi_{\K}(x-y)=x .
 \end{eqnarray}
The projection operator $\Pi_{\cal K}(x)$ is  known to be first-order directionally differentiable (see e.g. \cite[Lemma 2]{OS08}) and  the
 connection between its tangent cone and its  directional derivative has been  {given} (see \cite{JLZ15,YZ16}): for any $(x,y)\in \Omega$,
\begin{equation}
\label{directionalderivative}
 (d, w) \in T_\Omega(x,y) \Longleftrightarrow  \Pi'_{\cal K}(x-y;d-w)=d.
 \end{equation}
 Using this connection, it has been shown that the SOC complementarity set $\Omega$ is geometrically derivable and
  the exact formula for its tangent cone is given; see, e.g., \cite[Theorem 5.1]{YZ16}. Moreover, the coderivative of the projection operator $\Pi_{\cal K}$ allows us to characterize the various normal cones as in \cite[Proposition 2.1]{yzhou} and show that the SOC complementarity set is not only geometrically derivable but also directionally regular \cite[Theorem 6.1]{YZ16}.
  So far by using the first-order variational analysis, it has been revealed that although the SOC complementarity set is neither a convex set nor the union  of finitely many polyhedral convex sets, it enjoys  certain nice properties that a convex set or  the union  of finitely many polyhedral convex sets  has. In this paper, we continue to investigate the second-order variational  properties of the SOC complementarity cone.
Our main contributions are as follows:
 \begin{enumerate}
  \item[$\bullet$] We  derive the exact formula for the  second-order directional derivative of the projection operator over second-order cone. We {further} establish the connection between the second-order tangent set and the second-order directional derivative of the projection operator: for any $(x,y)\in \Omega$ and $(d,w) \in T_\Omega(x,y)$,
\begin{equation}
\label{sdirectionalderivative}(p,q) \in T^2_\Omega ((x,y);(d,w)) \Longleftrightarrow \Pi_{\cal K}''(x-y;d-w,p-q)=p.
\end{equation}
 \item[$\bullet$] We show that the SOC complementarity set is second-order directionally differentiable (see Definition \ref{BS}). Note that this
     nice property is not even enjoyed  by a convex set (see \cite[Example 3.31]{BS}).
 \item[$\bullet$]
Using the  characterization (\ref{sdirectionalderivative}) and the precise formula for the second-order directional derivative of the projection operation over the second-order cone, we derive the exact formula for the second-order tangent set of the SOC complementarity set. Compared with the usual vector complementarity set, our research shows that the task of establishing the formula of second-order tangent set to the second-order cone complementarity set, { which has nonpolyhedral and nonconvex structure}, is not trivial.
 \item[$\bullet$]{ Based on the exact formula of the second-order tangent set of $\Omega$, we develop the second-order optimality conditions for the mathematical program with second-order cone complementarity constraints (SOCMPCC).}
 \end{enumerate}


We organize our paper as follows. Section 2 contains the preliminaries. In Section 3, we calculate the second-order directional derivative of the projection operator over the second-order cone. Section 4 is devoted to the exact formula of the second-order tangent set to the SOC complementarity set. {The second-order optimality conditions of SOCMPCC are discussed in Section 5}.

 \section{Preliminaries}
 In this section, we clarify the notation and recall some background materials. First,
 we denote by $\Re_+$ and $\Re_{++}$ the set of nonnegative scalars and
 positive scalars respectively, i.e., $\Re_+:=\{\alpha|\ \alpha\geq 0\}$ and
 $\Re_{++}:=\{\alpha|\ \alpha>0\}$. For a set ${ C}$,  denote by int$C$, cl$C$, bd$C$, co$C$, $C^c$ its interior, closure, boundary,  convex hull, and its complement, respectively. For a closed set $C\subseteq \Re^n$, { let
 $C^\circ$ and $\sigma(\cdot|C)$
 stand for the polar cone and the support function of $C$, respectively, i.e.,
 $C^\circ=\{v|\, \langle v,w \rangle\leq 0,\  \forall w\in C \}$
 and $\sigma(z|C)=\sup\{\langle z,x \rangle|\, x\in C\}$ for $z\in \Re^n$.} Denote by ${\rm lin}C$ the largest subspace $L$ such that $C+L\subseteq C$.
 For a vector $x=(x_1,x_2)\in \Re\times \Re^{n-1}$, we denote $x^\circ$ the polar set of the set $\{x\}$ and $\hat{x}:=(x_1, -x_2)$, the reflection of vector $x$ on the  $x_1$ axis. For  a nonzero vector $x$, we denote by $\bar{x}:=x/\|x\|$.
 {Let  $ o(\lambda):\mathbb{R}_+\rightarrow \mathbb{R}^m$ stand for} a mapping with the property that $o(\lambda)/\lambda\rightarrow 0$ when $\lambda \downarrow 0$.
  For a  mapping $F:\mathbb{R}^n\to \mathbb{R}^m$ and   vectors $x, d \in \mathbb{R}^n$, we denote by $\nabla F(x)\in \mathbb{R}^{m\times n}$ the Jacobian  of $F$ at $x$,
 by $\nabla^2 F(x)$ the second-order derivative of $F$ at $x$, and by $\nabla^2 F(x)(d,d)$ the quadratic form corresponding to $\nabla^2 F(x)$.
 The directional derivative of $F$ at $x$ in direction $d$ is defined as
$$F'(x;d):=\lim_{{t\downarrow 0}}\frac{F(x+td) -F(x)}{t},
$$ provided that the above limit exists. If $F$ is directionally differentiable at $x$ in direction $d$, its  parabolic second-order directional derivative is defined as
$$F''(x;d, w):=\lim_{{t\downarrow 0}}\frac{F(x+td+\frac{1}{2}t^2w)-F(x)-tF'(x;d)}{\frac{1}{2}t^2},
$$ provided that the above limit exists. Moreover {if the following limit exists}
$$F''(x;d, w)=\lim_{{t\downarrow 0}\atop {w'\to w}}\frac{F(x+td+\frac{1}{2}t^2w')-F(x)-tF'(x;d)}{\frac{1}{2}t^2},
$$
then  $F$ is
 said to be parabolical second-order directionally differentiable at $x$ in the direction $d$ in the sense of Hadamard. In general, the concept of parabolical second-order directional differentiability in the Hadamard sense is stronger than that of  parabolical second-order directional differentiability. However, when $F$  is locally
Lipschitz at $x$,  these two concepts coincide.
It is known that if $F$ is  parabolical second-order directional differentiable in the Hadamard sense at $x$ along $d,w$, then
 \[
 F(x+td+\frac{1}{2}t^2w+o(t^2))=F(x)+tF'(x;d)+\frac{1}{2}t^2 F''(x;d,w)+o(t^2).
 \]


 \begin{definition}[Tangent Cones]\label{sotangent}
Let $S\subseteq \mathbb{R}^m$ and $x\in S$. {The regular/Clarke, inner and
 (Bouligand-Severi) tangent/contingent cone} to $S$ at $x$ are defined respectively as
 \begin{eqnarray*}
   \widehat{T}_S(x)&:=& { \liminf\limits_{{x' \stackrel{S}{\to}x} \atop {t\downarrow 0}}\frac{S-x'}{t}=\Big\{ d\in \Re^m\,\Big|\, \forall \, t_k\downarrow 0,\, {x_k \stackrel{S}{\to}x}, \; \exists d_k\to d \
  \ {\rm with}\ \ x_k+t_kd_k\in S  \Big\}},\\
 T^i_S(x)&:=&  \liminf\limits_{t\downarrow 0}\frac{S-x}{t}
 = \Big\{d\in\Re^m \, \Big| \, \forall \, t_k\downarrow 0,\; \exists d_k\to d \ \ {\rm with}
 \ \ x+t_k d_k\in S\Big\},\\
 T_S(x)&:=& \limsup\limits_{t\downarrow 0}\frac{S-x}{t}
 = \Big\{d\in\Re^m \, \Big| \, \exists \, t_k\downarrow 0,\;d_k\to d \ \ {\rm with}
 \ \ x+t_k d_k\in S\Big\}.
 \end{eqnarray*}
 The inner and
 outer second-order tangent sets to $S$ at $x$ in direction $d$ are defined
 respectively as
\begin{eqnarray*}
 T_S^{i,2}(x; d) &:=& \left \{ w \in \Re^m \, \bigg | \,
 \dist \left( x+t d + \frac{1}{2}t^2w, S \right) = o(t^2), \ \ t\geq 0 \right \}
,\\
 T_S^2(x; d) &:= &\left \{ w \in \Re^m \, \bigg | \, \exists \ t_n \downarrow 0 \
 {\rm such \ that}\ \dist \left( x+t_n d + \frac{1}{2}t_n^2w, S \right)
 = o(t_n^2) \right \}.
\end{eqnarray*}
   \end{definition}
While for a nonconvex set $S$, the contingent cone $T_S(x)$ may be nonconvex, it is known that the regular/Clarke tangent cone $\widehat{T}_S(x)$ is  always  closed and convex.
By definition, since the distance function of a convex set is convex, it is easy to see that the inner second-order  tangent set is always convex when the set $S$ is convex. On the other hand, the outer second-order  tangent set may be nonconvex even when the set $S$ is convex (see  \cite[Example 3.35]{BS}).
Note that
$T_S^{i,2}(x; d)\subseteq T_S^{2}(x; d)$ and the outer  second-order tangent set $T_S^{2}(x; d)$ needs not be a cone (it may be empty; see e.g. an example in \cite[page 592]{RW98}). If $T_S^{i,2}(x; d)=T_S^{2}(x; d)$, we simply call $T_S^{2}(x; d)$ the second-order tangent set to $S$ in direction $d$.

\begin{definition}{\cite[Definition 3.32]{BS}} \label{BS}A set $S$ is said to be second-order directionally differentiable at $x\in S$ in a direction $d\in T_S(x)$, if
 $T^i_S(x)=T_S(x)$ and $T^{i,2}_S(x;d)=T_S^2(x;d)$.
 \end{definition}
\begin{definition}[Normal Cones]  Let $S\subseteq \mathbb{R}^m$ and $x\in S$. The regular/Fr\'{e}chet, limiting/Mordukhovich, and Clarke normal cone of $S$ at ${x}$ are defined respectively as
  \begin{eqnarray*}
 \widehat{{N}}_S({x}) &:=& \Big \{v\in \mathbb{R}^m\,|\
 \langle v, x'-{x}\rangle \leq o(\|x'-x\|) \ \forall  x'\in S \Big\},\\
N_S (x)&:=&\limsup\limits_{x' \stackrel{S}{\to} x}\widehat{N}_{S}(x')= \Big \{\lim_{k\rightarrow \infty} v_k|\  v_k \in \widehat{N}_S(x_k),\ \   x_k \stackrel{S}{\to} x \Big\},\\
N^c_S (x)&:=& clcoN_S(x).
\end{eqnarray*}
\end{definition}
\begin{lemma}[Tangent-Normal Polarity](see \cite[Theorem 6.28]{RW98}, \cite{Clarke1983}) \label{polarity}
For a closed set $S\subseteq \mathbb{R}^m$ and $x\in S$, $\widehat{T}_S(x)=(N_S(x))^\circ=(N_S^c(x))^\circ, \widehat{N}_S(x)=(T_S(x))^\circ$, $(\widehat{T}_S(x))^\circ=N_S^c(x)$.
\end{lemma}

We recall some known results concerning the second-order cone ${\cal K}$ in $\Re^m$.
 The topological interior and the boundary of ${\cal K}$ are
\begin{eqnarray*}
{\rm int} {\cal K}= \{(x_1,x_2) \in \mathbb{R} \times \mathbb{R}^{m-1}|x_1>\|x_2\|\} \ \ {\rm and}\ \  {\rm bd} {\cal K}=\{(x_1,x_2) \in \mathbb{R} \times \mathbb{R}^{m-1}|x_1=\|x_2\|\},
\end{eqnarray*}
respectively.
Similar to the eigenvalue decomposition of a matrix, for any given vector $x:=(x_1,x_2)\in \Re\times\Re^{m-1}$, $x$ can be decomposed as
 (see e.g \cite{FLT02})
 \[
 x = \lambda_1(x) u_{x}^{(1)} + \lambda_2(x) u_{x}^{(2)},
 \]
 where $\lambda_i(x)$  and $u_{x}^{(i)}$  for $i=1,2$ are
 the spectral values and the associated spectral vectors of $x$
 respectively, given by
 \[
 \lambda_i(x):=x_1+(-1)^i\|x_2\|  \quad  {\rm and}\quad
 u^{(i)}_x:=\left\{\begin{array}{ll}\frac{1}{2}(1, (-1)^i \bar{x}_2) & {\rm if} \ \ x_2\neq 0,\\
\frac{1}{2}(1,(-1)^iw) & {\rm if}\ \ x_2=0,\end{array}\right.
\]
with $w$ being a fixed unit vector in $\mathbb{R}^{m-1}$.

\begin{lemma}(see e.g. \cite[Proposition 2.2]{YZ15})\label{lem2.1}
For any $x,y \in bd {\cal K}\backslash \{0\}$, the following equivalence holds:
$$x^Ty=0 \Longleftrightarrow y=k\hat{x} \mbox{ with } k=y_1/x_1>0  \Longleftrightarrow y=k\hat{x} \mbox{ with } k\in \mathbb{R}_{++}.$$
\end{lemma}

 For a given real-valued function $f:\Re\to \Re$, we define the SOC function
 $f^{\rm soc}: \Re^m \to \Re^m$ as
 \begin{equation}
 f^{\rm soc}(z) := f \left( \lambda_1(z) \right) u^{(1)}_z
 + f \left( \lambda_2(z) \right) u^{(2)}_z.\label{SOCfcn1}
 \end{equation}
For $z\in \Re^m$, let $\Pi_{\K}(z)$ be the metric projection of $z$ onto  $\K$. Then by \cite{FLT02}, it can be calculated as \begin{equation}
 \Pi_{\K}(z)=\lambda_1(z)_+u^{(1)}_z
 + \lambda_2(z)_+ u^{(2)}_z,\label{SOCfcn2}
 \end{equation}
 where $\alpha_+:=\max\{\alpha,0\}$ is the nonnegative part of the number  $\alpha\in \mathbb{R}$.
 Hence the projection operator $\Pi_{\cal K}(\cdot)$ is an SOC function corresponds to the plus function $f(\alpha):=\alpha_+$.


 \section{{Second-order directional derivative of the projection operator over the second-order cone}}

 As commented in the introduction, there exists a close relationship between the second-order tangent set of the SOC complementarity set and the second-order directional derivative of the projection operator $\Pi_{\cal K}$; see (\ref{sdirectionalderivative}). Therefore, to obtain the exact formula of the second-order tangent set, we need to  calculate  the second-order directional derivative of the projection operator $\Pi_{\cal K}$.  This task is done in this section, which is of independent interest.
For the convenience of notations, we sometime use $\Phi(x)$ instead of $\bar{x}$ to stand for ${x}/{\|x\|}$ as $x \neq 0$. It is easy to verify (see e.g. \cite[Theorem 3.1]{ZTC15}) that
 $\Phi$ is second-order continuously
 differentiable at $x\neq 0$ with
 \begin{eqnarray*}
  {\nabla} \Phi(x) &=& ( I-\bar{x}\bar{x}^T) /{\|x\|},\\
  {\nabla}^2\Phi(x)(w,w) &=& -2 \frac{\bar{x}^Tw}{\|x\|^2} w
 + w^T \left (\frac{ 3 \bar{x}\bar{x}^T -I }{\|x\|^3} \right) wx \\
 &=& -2\frac{\bar{x}^Tw}{\|x\|}\nabla\Phi(x)(w)-\frac{1}{\|x\|}w^T\nabla \Phi(x)w \bar{x},
  \end{eqnarray*}
 where $I$ is the identity matrix in $\Re^{m\times m}$.

 Since the second-order cone ${\cal K}$ is a special circular cone $\mathcal{L}_{\theta}$ defined by
 $$\mathcal{L}_{\theta}:=\{(x_1,x_2)\in \Re\times \Re^{n-1}| \cos\theta \|x\|\leq x_1\}$$
 with $\theta=45^\circ$, the SOC function $f^{\rm soc}$ is a special case of the circular cone function  $f^{\mathcal{L}_{\theta}}$  studied in  \cite{ZTC15} with $\theta= 45^{\circ}$.  The following result follows from \cite[Theorem 3.3]{ZTC15} immediately.

 \begin{lemma} \label{second-order diff}
 Suppose that $f:\Re\to \Re$. Then,  the SOC function $f^{\rm soc}$ is parabolic second-order
 directionally differentiable at $x$ in the Hadamard sense if and only if
 $f$ is parabolic second-order directionally differentiable at $\lambda_i(x)$
 in the Hadamard sense for $i=1,2$. Moreover,
 \begin{description}
 \item[(i)] if $x_2=0$ and $d_2=0$, then
            \[
            (f^{\rm soc})^{''}(x;d,w)
            = f^{''} \left( x_1;d_1,w_1-\|w_2\| \right) u_w^{(1)}
            + f^{''} \left( x_1;d_1,w_1+\|w_2\| \right) u_w^{(2)};
            \]
 \item[(ii)] if $x_2=0$ and $d_2 \neq 0$, then
            \begin{eqnarray*}
            (f^{\rm soc})^{''}(x;d,w)             &=& f^{''} \left( x_1;d_1-\|d_2\|,w_1-\bar{d}_2^Tw_2 \right) u_d^{(1)} \\
            & & + f^{''} \left( x_1;d_1+\|d_2\|,w_1+\bar{d}_2^Tw_2 \right) u_d^{(2)} \\
            & & + \frac{1}{2}\Big(f'(x_1;d_1+\|d_2\|)
                - f'(x_1;d_1-\|d_2\|) \Big)
           \left (
           \begin{matrix}
            0\\
        {\nabla}{\Phi}(d_2)w_2 \end{matrix}\right );
            \end{eqnarray*}
 \item[(iii)] if $x_2\neq 0$, then
            \begin{eqnarray*}
            \lefteqn{ (f^{\rm soc})^{''}(x;d,w)}\\
            &=& f^{''}\left( x_1-\|x_2\|;d_1-\bar{x}_2^Td_2,w_1
                - \big[ \bar{x}_2^T w_2+d_2^T {\nabla}\Phi(x_2) d_2\big] \right)u_x^{(1)} \\
            & & + f^{''} \left( x_1+\|x_2\|;d_1+\bar{x}_2^Td_2,w_1
                + \big[ \bar{x}_2^T w_2+d_2^T {\nabla}\Phi(x_2)d_2\big] \right) u_x^{(2)} \\
            & & + \Big(f'(x_1+\|x_2\|;d_1+\bar{x}_2^Td_2)
                - f'(x_1-\|x_2\|;d_1-\bar{x}_2^Td_2)\Big) \left ( \begin{matrix}
            0\\
        {\nabla}{\Phi}(x_2)d_2 \end{matrix}\right ) \\
            & & + \frac{1}{2}\Big(f(x_1+\|x_2\|)-f(x_1-\|x_2\|)\Big)\left ( \begin{matrix}
            0\\
      \nabla {\Phi}(x_2)w_2+ {\nabla}^2{\Phi}(x_2)(d_2,d_2) \end{matrix}\right ).
            \end{eqnarray*}
 \end{description}
 \end{lemma}

Since the projection operator $\Pi_{\cal K}(\cdot)$ is the SOC function corresponding to the plus function $f(\alpha):=\alpha_+$, we will need the second-order directional derivative of the plus function.

 \begin{lemma}(see e.g. \cite{ZZX13}) \label{Lem3.4} Let $f(\alpha):=\alpha_+$ for $\alpha \in \Re$. Then $f$ is parabolic second-order directionally differentiable at $x$
 in the Hadamard sense and
 \begin{equation}\nonumber
 f'(x;d)=\left\{
 \begin{array}{ll}
 d       & {\rm if} \ x>0, \\
 d_{+} & {\rm if} \ x=0, \\
 0       & {\rm if} \ x<0,
 \end{array}
 \right. \ \ \ {\rm and}\
 \quad f^{''}(x;d,w)=\left\{
 \begin{array}{ll}
 w       & {\rm if} \ x>0 \ {\rm or} \ x=0, d>0, \\
 0       & {\rm if} \ x<0 \ {\rm or} \ x=0, d<0, \\
 w_{+} & {\rm if} \ x=d=0.
 \end{array}
 \right.
 \end{equation}
\end{lemma}

Since in the formula of the second-order directional derivative of the projection operator, we will need the tangent cone and the second-order tangent set for the set ${\cal K}$ and its polar ${\cal K}^\circ$, for convenience we summarize their formulas in the following two lemmas.
 \begin{lemma}\cite[Lemma 25 and Lemma 27]{BR05}\label{Lem3.3}  For any $x\in \K$, one has\begin{eqnarray*}
 T_{\cal K}(x)
 &=& \left\{
 \begin{array}{lll}
 & \mathbb{R}^m \ &{\rm if } \  x\in {\rm int}{\cal K};  \\
 & {\cal K} \ &{\rm if}\ x=0; \\
 & \{ d\in \mathbb{R}^m| -d_1+\bar{x}_2^T d_2\leq 0\} \ \ & {\rm if}\ x\in {\rm bd}{\cal K}\backslash\{0\}.
 \end{array}
\right.
 \end{eqnarray*}For any $x\in \K$ and $d\in T_\K(x)$,
 \[
 T_{\K}^2(x;d)=\left\{
 \begin{array}{ll}
 \Re^m  \ & \ {\rm if } \  d\in {\rm int}T_\K(x); \\
  T_{\K}(d) &\ {\rm if }\  x=0; \\
 \{w|\ w_2^Tx_2-w_1x_1\leq d_1^2-\|d_2\|^2\} & \ {\rm if } \ x\in {\rm bd}\K\backslash \{0\} \mbox{ and } d\in {\rm bd}T_{\K}(x).
 \end{array}
 \right.
 \]
 \end{lemma}

 { Applying  \cite[Lemma 25 and Lemma 27]{BR05} to $\K^\circ=-\K$ yields the following result.}

 { \begin{lemma}\label{Lem3.5}
 For $x\in \K^\circ$, one has
 \[
 T_{\K^\circ}(x)=\left\{\begin{array}{ll}
 \Re^m         &  \ {\rm if } \ x\in {\rm int}\K^\circ;\\
  \K^\circ &  \ {\rm if } \ x=0;\\
 \{d\in \mathbb{R}^m| d_1+\bar{x}_2^Td_2\leq 0\}  &  \ {\rm if } \ x\in {\rm bd}\K^\circ\backslash\{0\}.
 \end{array}\right.
 \]
 For $x\in \K^\circ$ and $d\in T_{\K^\circ}(x)$, one has
  \[
 T^2_{\K^\circ}(x;d)=\left\{\begin{array}{ll}
 \Re^m &  \ {\rm if } \ d\in {\rm int}T_{\K^\circ}(x);\\
 T_{\K^\circ}(d)         &  \ {\rm if } \ x=0;\\
\{w| w_2^Tx_2-w_1x_1\leq d_1^2-\|d_2\|^2\}  &   \ {\rm if } \ x\in {\rm bd}\K^\circ\backslash \{0\} \mbox{ and } d\in {\rm bd}T_{\K^\circ}(x).
 \end{array}\right.
 \]
 \end{lemma}}

We are now ready to give  the second-order directional derivative of the projection operator.

 \begin{theorem} \label{projection-second-order diff}
  The projection operator $\Pi_{\K}$  is parabolic second-order
 directionally differentiable in the Hadamard sense. Moreover, for any $x,d,w\in \Re^m$,  the second-order directional derivative can be calculated as in the following six cases.
 \begin{description}
   \item[Case (i)] $x\in {\rm int} \K$.
  $\Pi_{\K}^{''}(x;d,w)=w$.
   \item[Case (ii)] $x\in {\rm int}\K^\circ$.
$\Pi_{\K}^{''}(x;d,w)=0$.
   \item [Case (iii)] $x=0$.
   \[
   \Pi_{\K}^{''}(x;d,w)=\left\{\begin{array}{ll}
   w \ \ \ \ \ \ \ \ \ \ \ \ &   \ {\rm if } \  d\in {\rm int}\K, \\
   0 &   \ {\rm if } \ d\in {\rm int}\K^\circ,\\
   \frac{1}{2}\left(\begin{matrix}
 w_1+\bar{d}_2^Tw_2\\
 \left[w_1-\frac{d_1}{\|d_2\|} {\bar{d}_2}^Tw_2\right]\bar{d}_2+\left[1+\frac{d_1}{\|d_2\|}\right]w_2
 \end{matrix}\right) &   \ {\rm if } \ d\in (\K\cup\K^\circ)^c, \\
  w  &    \ {\rm if } \ d\in {\rm bd}\K\backslash\{0\}, w\in T_{\K}(d),\\
  \frac{1}{2}\left(\begin{matrix}
  w_1+\bar{d}_2^Tw_2\\
 2w_2+(w_1-\bar{d}_2^Tw_2)\bar{d}_2\end{matrix}\right)  &   \ {\rm if } \  d\in {\rm bd}\K\backslash\{0\}, w\notin T_{\K}(d), \\
  0  &   \ {\rm if } \ d\in {\rm bd}\K^\circ\backslash \{0\},  w\in T_{\K^\circ}(d),\\
  \frac{1}{2}(w_1+\bar{d}_2^Tw_2)\left(\begin{matrix}
  1 \\
  \bar{d}_2
  \end{matrix}\right)  &  \ {\rm if } \  d\in {\rm bd}\K^\circ\backslash \{0\},  w\notin T_{\K^\circ}(d),\\
   \Pi_{\K}(w) &  \ {\rm if } \  d=0.
   \end{array} \right.
   \]
   \item[Case (iv)] $x\in {\rm bd}\K\backslash \{0\}$.
   \begin{eqnarray*}
\lefteqn{\Pi_{\K}^{''}(x;d,w)=}\\
 && \left\{\begin{array}{l}
   w  \qquad  \qquad \qquad  \qquad    \ {\rm if } \  d\in {\rm int}T_{\K}(x), \ \ \ \ \ \ \ \  \\
   w   \qquad  \qquad  \qquad  \qquad  \ {\rm if } \ d\in {\rm bd}T_{\K}(x), w\in T^2_{\K}(x;d), \\
  \frac{1}{2}\left(\begin{matrix}
  w_1+\bar{x}_2^Tw_2+\frac{\|d_2\|^2-d_1^2}{\|x_2\|}\\
  \left[w_1-\bar{x}_2^Tw_2-\frac{\|d_2\|^2-d_1^2}{\|x_2\|}\right]\bar{x}_2+2w_2
  \end{matrix}\right)\qquad  \ {\rm if } \  d\in {\rm bd} T_{\K}(x), w\notin T^2_{\K}(x;d), \\
  \frac{1}{2}\left(
  \begin{matrix}
  w_1+\bar{x}_2^Tw_2+\frac{\|d_2\|^2-(\bar{x}_2^Td_2)^2}{\|x_2\|}\\
  \left[w_1-\bar{x}_2^Tw_2-\frac{\|d_2\|^2-3(\bar{x}_2^Td_2)^2+2d_1\bar{x}_2^Td_2}{\|x_2\|}\right]\bar{x}_2
  +2w_2+2\frac{d_1-\bar{x}_2^Td_2}{\|x_2\|}d_2
  \end{matrix}
  \right) \quad  \ {\rm if } \  d\in T_{\cal K}(x)^c.
   \end{array}\right.
   \end{eqnarray*}

   \item[Case (v)] $x\in {\rm bd}\K^\circ\backslash \{0\}$.
   \begin{eqnarray*}
   \lefteqn{\Pi_{\K}^{''}(x;d,w)=}\\
   &&\left\{\begin{array}{l}
  0  \qquad  \qquad  \qquad  \qquad  \ {\rm if } \    d\in {\rm int}T_{\K^\circ}(x), \\
  0 \qquad  \qquad  \qquad  \qquad  \ {\rm if } \      d\in {\rm bd}T_{\K^\circ}(x), w\in T^2_{\K^\circ}(x;d), \\
  \frac{1}{2}\left(w_1+\bar{x}_2^Tw_2+\frac{\|d_2\|^2-d_1^2}{\|x_2\|}\right)
  \left(\begin{matrix}
  1 \\ \bar{x}_2
  \end{matrix}\right) \qquad   \ {\rm if } \  d\in {\rm bd} T_{\K^\circ}(x), w\notin T^2_{\K^\circ}(x;d),\\
  \frac{1}{2}\left(\begin{matrix}
  w_1+\bar{x}_2^Tw_2+\frac{\|d_2\|^2-(\bar{x}_2^Td_2)^2}{\|x_2\|}\\
  \left[w_1+\bar{x}_2^Tw_2+\frac{\|d_2\|^2-3(\bar{x}_2^Td_2)^2-2d_1\bar{x}_2^Td_2}{\|x_2\|}\right]\bar{x}_2+
  2\frac{d_1+\bar{x}_2^Td_2}{\|x_2\|}d_2
  \end{matrix}\right) \qquad \ {\rm if } \    d\in T_{\K^\circ}(x)^c.
   \end{array}\right.
   \end{eqnarray*}
   \item[Case (vi)] $x\in (\K\cup\K^\circ)^c$.
     \end{description}
   \begin{eqnarray*}
   &&\Pi_{\K}^{''}(x;d,w)=\\
   &&\frac{1}{2}\left(
   \begin{matrix}
   w_1+\bar{x}_2^Tw_2+\frac{\|d_2\|^2-(\bar{x}_2^Td_2)^2}{\|x_2\|}\\
   \left[w_1-\frac{x_1}{\|x_2\|}\bar{x}_2^Tw_2-\frac{x_1}{\|x_2\|^2}\big(\|d_2\|^2-3(\bar{x}_2^Td_2)^2\big)
   -2d_1\frac{\bar{x}_2^Td_2}{\|x_2\|}\right]\bar{x}_2+2\frac{\|x_2\|d_1-x_1\bar{x}_2^Td_2}{\|x_2\|^2}d_2
   +\left[1+\frac{x_1}{\|x_2\|}\right]w_2\end{matrix}
   \right).
   \end{eqnarray*}
  \end{theorem}

 \beginproof
   By  (\ref{SOCfcn1})-(\ref{SOCfcn2}), the projection operator $\Pi_{\K}$ is the SOC function $f^{soc}$ with
 $f(t):=t_+$. Applying Lemmas \ref{second-order diff} and \ref{Lem3.4}  will give the parabolic second-order
 directional differentiability of $\Pi_{\K}$ in the Hadamard sense  and a formula for  $\Pi^{''}_\K$. However in some cases the formula obtained will still involve the plus operator $(\cdot)_+$. In this theorem we aim at  obtaining the exact formula as proposed. For  some cases, e.g., in the cases $x\in {\rm int}\K$; $x\in {\rm int}\K^\circ$; $x=0,d\in {\rm int}\K$; $x=0,d\in {\rm int}\K^\circ$; $x=0,d=0$, we can prove the results  by directly using the definition of second-order directional derivative.
 In some other cases,  e.g., in the cases $x=0, d\in {\rm bd}\K\backslash \{0\}$; $x=0,d\in {\rm bd}{\rm \K^\circ}\backslash \{0\}$; $x \in {\rm bd}{\K}\backslash\{0\}, d\in {\rm bd}T_{\K}(x)$; $x\in {\rm bd}\K^\circ\backslash \{0\},  d\in {\rm bd} T_{\K^\circ}(x)$, we can further use the representation of tangent cones in Lemmas \ref{Lem3.3} and \ref{Lem3.5} to obtain the proposed exact formula. For simplicity, we only prove some of the  cases. The others can be obtained by following similar arguments.

 \noindent {\bf Case  }
 $x\in {\rm int}\K$. In this case $ \Pi_\K(x)=x$, $\Pi_\K'(x;d)=d$ and $\Pi_{\K}(x+td+\frac{1}{2}t^2w)=x+td+\frac{1}{2}t^2w$ for $t>0$ sufficiently small. Hence
 $$\Pi_{\K}^{''}(x;d,w):=\lim_{t\downarrow 0} \frac{ \Pi(x+td+\frac{1}{2}t^2 w)- \Pi_\K(x)-t\Pi_\K'(x;d)}{\frac{1}{2}t^2}=w. $$

 \noindent {\bf Case } $x=0$ and $d\in {\rm int}\K$. In this case $\Pi_\K(x)=0$ and $ \Pi_\K'(x;d)=d$. Note that  \[\Pi_{\K}(x+td+\frac{1}{2}t^2w)=\Pi_{\K}(td+\frac{1}{2}t^2w)=
 td+\frac{1}{2}t^2w,\]
 for $t>0$ sufficiently small. Hence  $\Pi_{\K}^{''}(x;d,w)=w.$

  \noindent {\bf Case } $x=0$ and $ d=0$. It is obvious that $\Pi_\K(0)=0, \Pi_\K'(0;0)=0$ and
  $ \Pi_{\K}(x+td+\frac{1}{2}t^2w)= \Pi_{\K}(\frac{1}{2}t^2w)=\frac{1}{2}t^2\Pi_{\K}(w)$. Hence $\Pi_{\K}^{''}(x;d,w)=\Pi_{\K}(w)$.

 \noindent {\bf Case} $x=0$ and $d\in {\rm bd}\K\backslash \{0\}$. Then $\Pi_\K(x)=0$ and $d_1=\|d_2\|\neq 0$. Directly applying Lemmas \ref{second-order diff}(ii) and  \ref{Lem3.4} yield
\begin{eqnarray}
 \Pi_{\K}^{''}(x;d,w)& =&
   \frac{1}{2}(w_1-\bar{d}_2^Tw_2)_+\left (\begin{matrix}1\\ -\bar{d}_2\end{matrix} \right ) +  \frac{1}{2}(w_1+\bar{d}_2^Tw_2)\left (\begin{matrix}1\\ \bar{d}_2\end{matrix} \right )  + \left (\begin{matrix}0\\( I-\bar{d}_2\bar{d}_2^T)w_2\end{matrix} \right). \label{eqn3-3-4-1}
 \end{eqnarray}
 Recall from Lemma \ref{Lem3.3} that $w\in T_{\K}(d)$ if and only if $w_1\geq \bar{d}_2^Tw_2$. It follows from (\ref{eqn3-3-4-1}) that
  $$\Pi_{\K}^{''}(x;d,w)=\left \{ \begin{array}{ll}
  w& \mbox{ if  } w\in T_{\K}(d),\\
    \frac{1}{2}\left(\begin{matrix}
  w_1+\bar{d}_2^Tw_2\\
 2w_2+(w_1-\bar{d}_2^Tw_2)\bar{d}_2\end{matrix}\right)
 &\mbox{ if } w\not \in T_{\K}(d).\end{array} \right. $$

   \noindent {\bf Case} $x \in {\rm bd}{\K}\backslash\{0\}$ and $d\in {\rm bd}T_{\K}(x)$. Then $ x_1=\|x_2\|\not =0$.
 and $-d_1+\bar{x}_2^T d_2=0$. Directly applying Lemmas \ref{second-order diff}(iii) and \ref{Lem3.4} yield
  \begin{eqnarray}
 \lefteqn{ \Pi_{\K}^{''}(x;d,w)}\nonumber\\
 &=& \frac{1}{2}\left( w_1- \left[\bar{x}_2^Tw_2+\frac{\|d_2\|^2-d_1^2}{\|x_2\|}\right] \right)_+\left(\begin{matrix}1\\ -\bar{x}_2\end{matrix}\right)  + \frac{1}{2}\left(\begin{matrix}
  w_1+\bar{x}_2^Tw_2+\frac{\|d_2\|^2-d_1^2}{\|x_2\|}\\
  \left(w_1-\bar{x}_2^Tw_2-\frac{\|d_2\|^2-d_1^2}{\|x_2\|}\right)\bar{x}_2+2w_2
  \end{matrix}\right). \nonumber\\
  \label{eqn3-4-2}
            \end{eqnarray}
  Recall from Lemma \ref{Lem3.3} that $w\in T^2_{\K}(x;d)$ if and only if $w_2^Tx_2-w_1x_1\leq d_1^2-\|d_2\|^2 $.
  Hence it follows from (\ref{eqn3-4-2}) that $\Pi_{\K}^{''}(x;d,w)=w$ if $w\in T^2_{\K}(x;d)$  and
  \begin{eqnarray*}
  \Pi_{\K}^{''}(x;d,w)= \frac{1}{2}\left(\begin{matrix}
  w_1+\bar{x}_2^Tw_2+\frac{\|d_2\|^2-d_1^2}{\|x_2\|}\\
  \left(w_1-\bar{x}_2^Tw_2-\frac{\|d_2\|^2-d_1^2}{\|x_2\|}\right)\bar{x}_2+2w_2
  \end{matrix}\right)
  \end{eqnarray*}
  if $w\notin T^2_{\K}(x;d)$.
 \endproof

  \section{Second-order tangent set for the SOC complementarity set}

%

This section is devoted to deriving the exact formula for the second-order tangent set to the SOC complementarity set.  To this end, we first  build its connection with the second-order directional derivative of the projection operator $\Pi_\K$, whose existence is guaranteed by virtue of Theorem \ref{projection-second-order diff}.
 \begin{proposition}\label{relationship}
For any $(x,y)\in \Omega$ and
 $(d,w)\in T_{\Omega}(x,y)$, one has
 \[
 T^{i,2}_{\Omega}\big((x,y);(d,w)\big)=T^2_{\Omega}\big((x,y);(d,w)\big)
 = \big \{(p,q) \, | \, \Pi^{''}_{\K}(x-y;d-w,p-q)=p  \big \}.
 \]
 \end{proposition}

  \beginproof
 Since $T^{i,2}_{\Omega}\big((x,y);(d,w)\big)\subseteq T^{2}_{\Omega}\big((x,y);(d,w)\big)$, it suffices to show
 \[
 T^{2}_{\Omega}\big((x,y);(d,w)\big) \subseteq \Upsilon\big((x,y);(d,w) \big) \subseteq T^{i,2}_{\Omega}\big((x,y);(d,w)\big),
 \]
where  $
 \Upsilon\big((x,y);(d,w) \big):=\big \{(p,q) \, | \, \Pi^{''}_{\K}(x-y;d-w,p-q)=p  \big \}.
$
 Let $(p,q)\in T^{2}_{\Omega}\big((x,y);(d,w)\big)$. Then by definition, there exist $t_n \downarrow 0$, $(\alpha(t_n),\beta(t_n))=o(t_n^2)$ such that
  $(x,y)+t_n(d,w) + \frac{1}{2}t_n^2 (p,q)
     + (\alpha(t_n),\beta(t_n))\in \Omega.$
 By the equivalence in (\ref{Omeganew}), it follows that
 \begin{eqnarray*}
\lefteqn{ \Pi_{\K} \Big(x-y+t_n (d-w)+\frac{1}{2}t_n^2 (p-q)+\alpha(t_n)-\beta(t_n)\Big)} \\
 &=& x+t_nd + \frac{1}{2}t_n^2p + \alpha(t_n) \\
 &=& \Pi_{\K}(x-y) + t_n \Pi'_{\cal K}(x-y;d-w) + \frac{1}{2}t_n^2p + \alpha(t_n),
 \end{eqnarray*}
 where the last equality follows from the equivalence in (\ref{Omeganew}) and
(\ref{directionalderivative}).
 Hence, $\Pi^{''}_{\K}(x-y;d-w,p-q)=p$, i.e., $(p,q)\in \Upsilon\big((x,y);(d,w) \big)$.

 Now, take $(p,q)\in \Upsilon\big((x,y);(d,w) \big)$, i.e., $ \Pi^{''}_{\K}(x-y;d-w,p-q)=p$. For $t>0$, define
 \[
 r(t):=\Pi_{\K}(x-y+t(d-w)+\frac{1}{2}t^2(p-q))-
 \Pi_{\K}(x-y)-t\Pi'_{\K}(x-y;d-w)-\frac{1}{2}t^2\Pi^{''}_{\K}(x-y;d-w,p-q).
 \]
 Then $r(t)=o(t^2)$ according to the second-order directional differentiability of $\Pi_{\K}$ by Theorem \ref{projection-second-order diff}. Note that
 \begin{eqnarray*}
 &&\Pi_{\K}\big(x-y+t(d-w)+\frac{1}{2}t^2p+r(t)-\frac{1}{2}t^2q-r(t)
 \big)\\
 &=&\Pi_{\K}(x-y+t(d-w)+\frac{1}{2}t^2(p-q))\\
 &=& \Pi_{\K}(x-y)+t\Pi'_{\K}(x-y;d-w)+\frac{1}{2}t^2
 \Pi_{\K}^{''}(x-y;d-w,p-q)+r(t)\\
 &=& x+td+\frac{1}{2}t^2p+r(t),
 \end{eqnarray*}
   where the last equality follows from the equivalence in (\ref{Omeganew}) and
(\ref{directionalderivative}).
  This together with  equivalence  (\ref{Omeganew})  yields
that
  \[
 \left(x+td+\frac{1}{2}t^2p+r(t), y + tw +
 \frac{1}{2}t^2q+r(t) \right)\in \Omega.
 \]
 It means $(p,q)\in  T^{i,2}_{\Omega}\big((x,y);(d,w)\big)$. The proof is complete.
 \endproof

  \begin{remark} The proof of equivalence (\ref{sdirectionalderivative}) in Proposition \ref{relationship}  is very similar to that of equivalence (\ref{directionalderivative}) as in  \cite[Proposition 5.2]{YZ16}. Note that although the equivalence (\ref{directionalderivative}) was shown in \cite[Proposition 3.1]{JLZ15}, the proof in \cite[Proposition 5.2]{YZ16} is much more concise without going over each possible cases as in \cite[Proposition 3.1]{JLZ15}. Moreover from the proof of \cite[Proposition 5.2]{YZ16}, one can {see} that the equivalence (\ref{directionalderivative}) holds for any general convex cone $\K$  as long as the  projection operator $\Pi_{\K}$ satisfies the Lipschitz continuity and directional differentiability.
 Similarly from the proof  of Proposition \ref{relationship}, we can see that  equivalence (\ref{sdirectionalderivative}) in Proposition \ref{relationship} holds for any general convex cone $\K$  whenever the projection operator $\Pi_{\K}$ satisfies the Lipschitz continuity and parabolic second-order directional differentiability in the Hadamard sense.
 \end{remark}

 The above result tells us that for characterizing the structure of the
 second-order tangent set to $\Omega$, we need to study the expression
 of the second-order directional derivative of the projection operator
 $\Pi_{\K}$, which has been obtained in Theorem \ref{projection-second-order diff}. With these preparations, the explicit expression of the second-order
 tangent set to $\Omega$ is given below. For convenience, we recall the formula for the tangent cone first.

 \begin{lemma}\label{tangentcone}\cite[Theorem 5.1]{YZ16}
For any $(x,y)\in \Omega$,
 \begin{eqnarray*}
 \lefteqn{T^{i}_\Omega(x,y)=T_\Omega(x,y)}\\
 &=& \left\{ (d,w) \left|
 \begin{array}{lll}
 & d\in \mathbb{R}^m,\ w=0, \ &{\rm if}\ x\in {\rm int}{\cal K}, \ y=0; \\
  & d=0,\ w\in \mathbb{R}^m, \ &{\rm if} \ x=0, \ y\in {\rm int}{\cal K};  \\
 & x_1\hat{w}-y_1d\in \mathbb{R} x, \ \ d\perp y,\ w\perp x, \ \ & {\rm if}
   \ x,y \in {\rm bd}{\cal K}\backslash\{0\}; \\
 & d\in T_{\cal K}(x),\ w=0 \ {\rm or}\ d\perp \hat{x}\ , w\in \mathbb{R}_{+}\hat{x}, \ \
 & {\rm if}\ x\in {\rm bd}{\cal K}\backslash \{0\},\ y=0; \\
 & d=0,\ w\in T_{\cal K}(y) \ {\rm or}\ d\in \mathbb{R}_{+}\hat{y}, \ w\perp \hat{y}, \ \
 & {\rm if}\ x=0,\ y\in {\rm bd}{\cal K}\backslash \{0\}; \\
 & d\in {\cal K},\ w\in {\cal K}, \ d\perp w, \ \ & {\rm if}\ x=0, \ y=0.
 \end{array}
 \right.\right\}.
 \end{eqnarray*}
 \end{lemma}

 According to Proposition \ref{relationship} and Lemma \ref{tangentcone}, we obtain the following result.
 \begin{theorem}\label{second-order-directionally}
 The set $\Omega$ is second-order directionally differentiable at every $(x,y)\in \Omega$ in every direction $(d,w)\in T_{\Omega}(x,y)$.
 \end{theorem}

  \begin{remark}
 It is well-known that for a convex set, the tangent cone and inner tangent cone coincide, but the inner and outer second-order tangent sets can be different; see \cite[Example 3.31]{BS}. Here we show that SOC complementarity set $\Omega$, although it is nonconvex, is second-order directionally differentiable, i.e., the tangent cone and inner tangent cone coincide, and the inner and outer  second-order tangent sets coincide as well.
 \end{remark}

 The inner and outer second-order tangent set to product sets have been studied in \cite[Page 168]{BS}. Particularly, for $C:=C_1\times\cdots\times C_m$ with $C_i\in \Re^{n_i}$, at certain $x=(x_1,\dots,x_m)$ with $x_i\in C_i$, according to \cite{BS},
 \[
  T^{i,2}_{C}(x;d)=T^{i,2}_{C_1}(x_1;d_1)\times\dots T^{i,2}_{C_m}(x_m;d_m)
  \]
  and
  \begin{equation}\label{R-product-1}
  T^{2}_{C}(x;d)\subset T^{2}_{C_1}(x_1;d_1)\times\dots T^{2}_{C_m}(x_m;d_m).
  \end{equation}
  If all except at most one of $C_i$ are second-order directional differentiable, then the equality holds in (\ref{R-product-1}).
Noting that second-order cone complementarity set is second-order directional differentiable, Theorem  \ref{second-order-directionally} can be then extended to the Cartesian product of finitely many second-order cone complementarity sets.

   \begin{corollary}
   Suppose that $\Omega_1,\cdots,\Omega_l$ are all SOC complementarity sets. Then the Cartesian product $\Omega:=\Omega_1\times \Omega_2\times \cdots \times \Omega_l$ is second-order directionally differentiable at every $(x,y)\in \Omega$ in every direction $(d,w)\in T_{\Omega}(x,y)$ and
 $$T^2_{\Omega}((x,y);(d,w))=T^2_{\Omega_1}((x_1,y_1);(d_1,w_1))\times \cdots \times T^2_{\Omega_1}((x_l,y_l);(d_l,w_l)).$$
 \end{corollary}

 \beginproof
 Since $(d,w)\in T_{\Omega}(x,y)=T_{\Omega_1} (x_1,y_1)\times \dots \times T_{\Omega_l} (x_l,y_l)$,  $(d_i,w_i)\in T_{\Omega_i}(x_i,y_i)$ for $i=1,\dots,l$. Take $(p,q)\in T^2_{\Omega}((x,y);(d,w))$. Hence
 \begin{eqnarray*}
 T^2_{\Omega}((x,y);(d,w)) &\subseteq & T^2_{\Omega_1}((x_1,y_1);(d_1,w_1))\times \cdots \times T^2_{\Omega_1}((x_l,y_l);(d_l,w_l)) \\
 & =& T^{i,2}_{\Omega_1}((x_1,y_1);(d_1,w_1))\times \cdots \times T^{i,2}_{\Omega_l}((x_l,y_l);(d_l,w_l))\\
 &= & T^{i,2}_{\Omega}((x,y);(d,w)),
 \end{eqnarray*}
  where the first inclusion and the second equation follows from \cite[Page 168]{BS}, and the first equation comes from Theorem \ref{second-order-directionally}.
   \endproof

 \begin{theorem} \label{formula-regular normal cone-1}
For any $(x,y)\in \Omega$ and
 $(d,w)\in T_\Omega(x,y)$, the formula of the second-order tangent set for the SOC complementarity set can be described as in the following six cases.

 \begin{description}
 \item[Case (i)]  $ x\in {\rm int} \K\ {\rm and } \ y=0$.
 $ T^2_\Omega\big((x,y);(d,w)\big)=\Re^m \times \{0\}. $
 \item[Case (ii)] $ x=0 \ {\rm and }\  y \in {\rm int} \K$.
 $ T^2_\Omega\big((x,y);(d,w)\big)=\{0\}\times \Re^m  .$
\item[Case (iii)] $x,y\in {\rm bd}\K\backslash \{0\}$.
 \begin{eqnarray*}
\lefteqn{ T^2_\Omega\big((x,y);(d,w)\big)}\\
&& = \left\{(p,q)\left|
 \begin{array}{ll}
  p\in {\rm bd} T^2_{\K}(x;d), \ \ q\in {\rm bd}T^2_{\K}(y;w), \\
  (x_1w_1-y_1d_1)\left(\frac{w_2-w_1\bar{y}_2}{y_1}-\frac{d_2-d_1\bar{x}_2}{x_1}\right)-
 p_1y_2-q_1x_2=x_1q_2+y_1p_2
 \end{array} \right.\right\}.
 \end{eqnarray*}

\item[Case (iv)] $x\in {\rm bd}\K\backslash \{0\}$ and $y=0$.
 \begin{eqnarray*}
 \lefteqn{T^2_\Omega\big((x,y);(d,w)\big)} \\
 &=&\left\{(p,q)\left|
 \begin{array}{lll}
 q=0,  \ & {\rm if} \ d\in {\rm int} T_{\K}(x), \ w=0; \\
 p\in T_{\K}^2(x;d),\ q=0, \ {\rm or} \ \
 p\in {\rm bd} T^2_{\K}(x;d), \ q\in \Re_+\hat{x}\ & {\rm if}\ d\in {\rm bd}T_{\K}(x),\ w=0; \\
 p \in {\rm bd}T^2_{\K}(x;d),\ -q_1\bar{x}_2-2\frac{w_1d_2}{\|x_2\|}
 -2 \frac{d_1w_2}{\|x_2\|}=q_2, \ \ & {\rm if}\ d\perp \hat{x}, \ w\in \Re_{++}\hat{x}.
 \end{array}
 \right.\right\}.
 \end{eqnarray*}
\item[Case (v)] $x=0$ and $y\in {\rm bd}\K\backslash \{0\}$.
 \begin{eqnarray*}
 \lefteqn{ T^2_\Omega\big((x,y);(d,w)\big)}\\
 &=& \left\{(p,q)\left|
 \begin{array}{lll}
 p=0,  \ & {\rm if} \ d=0,\ w\in {\rm int} T_{\K}(y); \\
 p=0,\ q\in T_{\K}^2(y;w), \ {\rm or} \  p\in \Re_+\hat{y},
 \ \ q\in {\rm bd} T^2_{\K}(y;w)\ & {\rm if}\ d=0,\ w\in {\rm bd}T_{\K}(y); \\
 q\in {\rm bd}T^2_{\K}(y;w),\ -p_1\bar{y}_2-2\frac{w_1d_2}{\|y_2\|}
 - 2\frac{d_1w_2}{\|y_2\|}=p_2, \ \ & {\rm if}\ d\in \Re_{++}\hat{y}, \ w\perp \hat{y}.
 \end{array}
 \right.\right\}.
 \end{eqnarray*}

\item[Case (vi)]   $ x=y=0$.
 $ T^2_\Omega\big((x,y);(d,w)\big)= T_\Omega(d,w) .$
 \end{description}
   \end{theorem}

 \beginproof
  By Proposition \ref{relationship}, to describe an element $(p,q)\in  T_\Omega^2((x,y);(d,w))$, it suffices to
 describe an element $(p,q)$ satisfying $\Pi_\K^{''}(x-y;d-w,p-q)=p$.
 For simplicity, we  denote by
$
 z:=x-y,\xi:=d-w $ and $  \eta:=p-q.$

 \medskip
 \noindent
{\bf Case (i)}  $x\in {\rm int}\K$ and $y=0$. Since $z=x-y\in {\rm int}\K$,  by Theorem \ref{projection-second-order diff}(i), we have $\Pi_{\K}^{''}(x-y;d-w,p-q)=p-q$.
It follows that
$$\Pi_\K^{''}(x-y;d-w,p-q)=p\Longleftrightarrow q=0.$$
 Hence $ T^2_\Omega\big((x,y);(d,w)\big)=\Re^m \times \{0\}. $

 \medskip
 \noindent
 {\bf Case (ii)} $x=0$ and $y\in {\rm int}\K$.  Since $z=x-y\in -{\rm int}\K$, by Theorem \ref{projection-second-order diff}(ii), we know
$
 \Pi^{''}_{\K}(z;d-w,p-q)=0.$
It follows that
$$\Pi_\K^{''}(x-y;d-w,p-q)=p\Longleftrightarrow p=0.$$
Hence $ T^2_\Omega\big((x,y);(d,w)\big)=\{0\} \times \Re^m. $

 \medskip
 \noindent
{\bf Case (iii)} $x,y\in {\rm bd}\K\backslash \{0\}$ and $x^Ty=0$. In this case $x_1=\|x_2\|\not =0$ and  by Lemma \ref{lem2.1},
 \begin{equation}
 z=x-y=(x_1,x_2)-k(x_1,-x_2)=((1-k)x_1,(1+k)x_2), \ \ k=y_1/x_1. \label{case3}\end{equation}
 This yields $z_1+\|z_2\|=2x_1>0$ and $z_1-\|z_2\|=-2kx_1<0$, i.e., $z\in (\K\cup \K^\circ)^c$.
Then by
 Theorem \ref{projection-second-order diff}(vi), $\Pi_{\K}^{''}(z;\xi,\eta)=p$ where $p=(p_1,p_2)\in \mathbb{R}\times \mathbb{R}^{m-1}$ if and only if
 \begin{eqnarray}
p_1 &=& \frac{1}{2}\left(\eta_1+\bar{z}_2^T\eta_2
 +\frac{\|\xi_2\|^2-(\bar{z}_2^T\xi_2)^2}{\|z_2\|}\right), \label{formula-case-3-1}
 \\
  p_2&=&\frac{1}{2}\left(\eta_1-\frac{z_1}{\|z_2\|}\bar{z}_2^T\eta_2
  -\frac{z_1}{\|z_2\|^2}
  \Big[\|\xi_2\|^2-3(\bar{z}_2^T\xi_2)^2\Big]
   -2\xi_1\frac{\bar{z}_2^T\xi_2}{\|z_2\|}\right)\bar{z}_2\nonumber\\
   && +\frac{\|z_2\|\xi_1-z_1\bar{z}_2^T\xi_2}{\|z_2\|^2}\xi_2
   +\frac{1}{2}\left(1+\frac{z_1}{\|z_2\|}\right)\eta_2.\label{formula-case-3-2}
 \end{eqnarray}

  We now try to derive an equivalent expression for  (\ref{formula-case-3-1}) and (\ref{formula-case-3-2}). Since $(d,w)\in T_{\Omega}(x,y)$, according to Lemma \ref{tangentcone},
  $ x \perp w$, $y\perp d$ and there exists $\beta\in \Re$ such  that $x_1\hat{w}-y_1d=\beta x$, from which and $x_1=\|x_2\|\not =0$  we have
 \begin{equation} \label{3-6}
w_1=kd_1+\beta, \quad w_2=-kd_2-\beta\bar{x}_2,
 \end{equation}
   and
 \begin{equation}
 \bar{x}_2^Tw_2=-w_1, \quad \bar{x}_2^Td_2=d_1.\label{eqn3-12}
 \end{equation}
  Note that $\bar{z}_2=\bar{x}_2$ by (\ref{case3}). Hence it follows from  (\ref{3-6}) and (\ref{eqn3-12}) that
  \begin{eqnarray}
 && \bar{z}_2^T\xi_2 = \bar{x}_2^T(d_2-w_2)=d_1+w_1 = (1+k)d_1 + \beta,\label{eqn3-13}\\
 &&\|\xi_2\|^2=\|d_2-w_2\|^2
 = \|(1+k)d_2+\beta \bar{x}_2\|^2
 = (k+1)^2\|d_2\|^2+2\beta(k+1)d_1 + \beta^2,\nonumber \\
&& \xi_1 = d_1-w_1 = (1-k)d_1-\beta. \label{eqn3-19}
   \end{eqnarray}
 Hence (\ref{formula-case-3-1}) can be rewritten as
\begin{equation} \label{3-3old}
 p_1=-q_1+\bar{x}_2^T(p_2-q_2)+\frac{x_1+y_1}{x_1^2}\Big(\|d_2\|^2-d^2_1\Big). \end{equation}
 The term in front of $\bar{z}_2$ in (\ref{formula-case-3-2}) becomes
 \begin{eqnarray*}
 &&\frac{1}{2}\left(\eta_1-\frac{z_1}{\|z_2\|}\bar{z}_2^T\eta_2-\frac{z_1}{\|z_2\|^2}
  \Big[\|\xi_2\|^2-3(\bar{z}_2^T\xi_2)^2\Big]
   -2\xi_1\frac{\bar{z}_2^T\xi_2}{\|z_2\|}\right)\\
 &= & \frac{1}{2}\left(\eta_1+\bar{z}_2^T\eta_2
     + \frac{\|\xi_2\|^2-(\bar{z}_2^T\xi_2)^2}{\|z_2\|}\right)
     -\frac{z_1+\|z_2\|}{2\|z_2\|}\bar{z}_2^T\eta_2\\
     && +\frac{z_1+\|z_2\|}{2\|z_2\|^2}\Big[(\bar{z}_2^T\xi_2)^2-\|\xi_2\|^2\Big] +\left[\frac{z_1}{\|z_2\|^2}\big(\bar{z}_2^T\xi_2\big)^2
     -\frac{\xi_1}{\|z_2\|}\big(\bar{z}_2^T\xi_2\big)\right] \\
   &=& \frac{y_1p_1-x_1q_1}{x_1+y_1}+\left[\frac{x_1-y_1}{(x_1+y_1)^2}(d_1+w_1)^2
     -\frac{1}{x_1+y_1}(d_1^2-w_1^2) \right] \\
 &=& \frac{y_1p_1-x_1q_1}{x_1+y_1}
     +2\frac{x_1w_1^2+x_1d_1w_1-y_1d_1^2-y_1d_1w_1}{(x_1+y_1)^2},
 \end{eqnarray*}
 where the second equality  uses (\ref{case3}), (\ref{formula-case-3-1}), and  (\ref{eqn3-13})-(\ref{3-3old}).
 It follows from (\ref{eqn3-13}) and (\ref{eqn3-19}) that  the term in front of $\xi_2$ in (\ref{formula-case-3-2}) is
 \[
 \frac{\xi_1}{\|z_2\|}- \frac{z_1}{\|z_2\|}\frac{\bar{z}_2^T\xi_2}{\|z_2\|}
 = \frac{1}{x_1+y_1}\big(d_1-w_1\big)-\frac{x_1-y_1}{(x_1+y_1)^2}\big(d_1+w_1\big)
 = 2\frac{y_1d_1-x_1w_1}{(x_1+y_1)^2}.
 \]
 The term in front of  $\eta_2$ in (\ref{formula-case-3-2}) is
 $
 1/2\big(1+(z_1/\|z_2\|)\big)=x_1/(x_1+y_1).
 $
 Hence (\ref{formula-case-3-2}) can be rewritten as
 \begin{eqnarray}
p_2 &= & \left(\frac{y_1p_1-x_1q_1}{x_1+y_1}
     +2\frac{x_1w_1^2+x_1d_1w_1-y_1d_1^2-y_1d_1w_1}{(x_1+y_1)^2}\right)\bar{x}_2+2\frac{y_1d_1-x_1w_1}{(x_1+y_1)^2} (d_2-w_2)\nonumber\\
     &&
    +\frac{x_1}{x_1+y_1}(p_2-q_2).\label{eqn4-3-2}
 \end{eqnarray}
 Further notice that
 $
 (y_1p_1-x_1q_1)\bar{x}_2=-p_1y_2-q_1x_2
 $
 and
 \begin{eqnarray*}
\lefteqn{2 \frac{y_1d_1-x_1w_1}{(x_1+y_1)^2}(d_2-w_2)
     +2 \frac{x_1w_1^2+x_1d_1w_1-y_1d_1^2-y_1d_1w_1}{(x_1+y_1)^2}\bar{x}_2 } \nonumber\\
  &=& 2 \frac{-x_1\beta}{(x_1+y_1)^2}[(1+k)d_2+\beta\bar{x}_2]
     +2 \frac{x_1\beta^2+y_1d_1\beta+x_1d_1\beta}{(x_1+y_1)^2}\bar{x}_2 \nonumber \\
  &=& 2\frac{\beta}{x_1+y_1}\big(-d_2+d_1\bar{x}_2\big)
  = 2\frac{x_1w_1-y_1d_1}{x_1(x_1+y_1)}\big(-d_2+d_1\bar{x}_2\big),
 \end{eqnarray*}
 where the first and third equations use (\ref{3-6}) and  the fact $y_1=kx_1$. Hence (\ref{eqn4-3-2}) can be rewritten as
  \begin{eqnarray}
  y_1p_2+x_1q_2&=&2 \frac{x_1w_1-y_1d_1}{x_1}(-d_2+d_1\bar{x}_2)-p_1y_2-q_1x_2 \nonumber\\
  &=& (x_1w_1-y_1d_1)\left(\frac{w_2-w_1\bar{y}_2}{y_1}+\frac{-d_2+d_1\bar{x}_2}{x_1}\right)-p_1y_2-q_1x_2, \label{eqn4-3-1}
  \end{eqnarray}
  where the second step comes from the fact
    $ (-d_2+d_1\bar{x}_2)/x_1=(w_2-w_1\bar{y}_2)/y_1$ due to (\ref{3-6}).
  Hence  (\ref{formula-case-3-1}) and (\ref{formula-case-3-2}) is equivalent to
 (\ref{3-3old}) and (\ref{eqn4-3-1}).

  Now, multiplying (\ref{eqn4-3-1}) by $\bar{x}_2^T$ and using (\ref{eqn3-12})  yields
 \begin{equation} \label{eqn4-3-4}
x_1\bar{x}_2^Tq_2+y_1\bar{x}_2^Tp_2= y_1p_1-x_1q_1.
 \end{equation}
 Hence it follows from (\ref{3-3old}) and (\ref{eqn4-3-4}) that
 \[
  p_1= \left(1+\frac{y_1}{x_1}\right)\bar{x}_2^Tp_2-\frac{y_1}{x_1}p_1+\frac{x_1+y_1}{x_1^2}\big(\|d_2\|^2-d^2_1\big),
  \]
   i.e.,
   \begin{equation}\label{eqn4-3-5} p_1=\bar{x}_2^Tp_2+\frac{1}{x_1}\big(\|d_2\|^2-d_1^2\big) \Longleftrightarrow p\in {\rm bd}T_\K^2(x; d).
 \end{equation}
 Since $d_1=(w_1-\beta)/k$, $d_2=-(w_2+\beta\bar{x}_2)/k$
 by (\ref{3-6}), and $w_2^T\bar{x}_2=-w_1$ by (\ref{eqn3-12}), we see
 \begin{equation}\label{case-3-a}
 \frac{\|d_2\|^2-d_1^2}{x_1^2}
 =\frac{\|w_2\|^2-w_1^2}{y_1^2}.
 \end{equation}
  Similarly, it follows from (\ref{3-3old}), (\ref{eqn4-3-4}), (\ref{case-3-a}), and $\bar{x}_2=-\bar{y}_2$ that
 \begin{eqnarray*}
 q_1=
 -\frac{x_1}{y_1}q_1+\left(1+\frac{x_1}{y_1}\right)\bar{y}_2^Tq_2+\frac{x_1+y_1}{y_1^2}
 \big(\|w_2\|^2-w_1^2\big),
 \end{eqnarray*}
 i.e.,
 \begin{equation}\label{eqn4-3-6}
 \frac{1}{y_1}\big(\|w_2\|^2-w_1^2\big)=q_1-\bar{y}_2^Tq_2 \Longleftrightarrow q\in {\rm bd} T_\K^2(y; w).
 \end{equation}
Hence along the line
 \[
 \left\{\begin{array}{l} (\ref{formula-case-3-1})\\
 (\ref{formula-case-3-2}) \end{array}\right. \Longleftrightarrow \left\{\begin{array}{l} (\ref{3-3old})\\
(\ref{eqn4-3-1}) \end{array}\right.\Longleftrightarrow  \left\{\begin{array}{l} (\ref{eqn4-3-5}), (\ref{eqn4-3-6})\\
(\ref{eqn4-3-1}) \end{array}\right.
 \]
 the desired result follows.

 \medskip
 \noindent
{\bf Case (iv)}  $x\in {\rm bd}\K\backslash \{0\}$ and $y=0$. In this case
 $z=x-y=x\in {\rm bd}\K\backslash \{0\}$.

 \medskip

$ {\bf (iv)}$-1. $d\in {\rm int} T_{\K}(x)$ and $w=0$.
 Then $\xi=d-w=d\in {\rm int}T_{\K}(x)$. Hence $\Pi_{\K}^{''}(z;\xi,\eta)=\eta$  by Theorem  \ref{projection-second-order diff}(iv). It follows that $\Pi_{\K}^{''}(x-y;d-w,p-q)=p$ if and only if $q=0$.

 \medskip
 $ {\bf (iv)}$-2. $d\in {\rm bd} T_{\K}(x)$ and $w=0$. Then
 $\xi=d\in {\rm bd} T_{\K}(x)$.
Hence $\Pi^{''}_{\K}(z;\xi,\eta)=\Pi_{\K}^{''}(x;d,\eta)$ and by Proposition \ref{second-order diff}(iv)
 \begin{eqnarray*}
\lefteqn{\Pi_{\K}^{''}(x;d,\eta)}\\
 &=& \left\{\begin{array}{ll}
   \eta \ \ \ \ \ \ \ \ \qquad  &  \mbox{if }  \eta \in T^2_{\K}(x;d), \\
  \frac{1}{2}\left(\begin{matrix}
  \eta_1+\bar{x}_2^T\eta_2+\frac{\|d_2\|^2-d_1^2}{\|x_2\|}\\
  \left(\eta_1-\bar{x}_2^T\eta_2-\frac{\|d_2\|^2-d_1^2}{\|x_2\|}\right)\bar{x}_2+2\eta_2
  \end{matrix}\right) & \mbox{if }  \eta \notin T^2_{\K}(x;d).
   \end{array}\right.
   \end{eqnarray*}
Note that
$\eta \in T^2_{\K}(x;d)\Longleftrightarrow\eta_2^Tx_2-\eta_1x_1\leq d_1^2-\|d_2\|^2$ by Lemma \ref{Lem3.3}.
 Hence $\Pi_{\K}^{''}(x;d,p-q)=p$ if and only if either $p-q\in T^2_{\K}(x;d)$ and $q=0$ or the following system holds
     \begin{eqnarray}\label{Case-4-Part-2}
  \left\{ \begin{array}{ll}
     \eta_1-\bar{x}_2^T\eta_2-\frac{\|d_2\|^2-d_1^2}{\|x_2\|} < 0, \\
     \frac{1}{2}\left(\eta_1+\left[\bar{x}_2^T\eta_2
     + \frac{\|d_2\|^2-d_1^2}{\|x_2\|} \right] \right)= p_1, \\
     \frac{1}{2}\left(\eta_1-\left[\bar{x}_2^T\eta_2
     +\frac{\|d_2\|^2-d_1^2}{\|x_2\|}\right]\right)\bar{x}_2
     +\eta_2=p_2. \end{array}\right.
 \end{eqnarray}
We now further simplify the  system (\ref{Case-4-Part-2}).
 \begin{eqnarray*}
 (\ref{Case-4-Part-2}) & \Longleftrightarrow&   \left\{ \begin{array}{ll}
    \eta_1-\left[\bar{x}_2^T\eta_2+\frac{\|d_2\|^2-d_1^2}{\|x_2\|}\right]< 0 ,\\ \bar{x}_2^T\eta_2+\frac{\|d_2\|^2-d_1^2}{\|x_2\|} = p_1 + q_1, \\
     \frac{1}{2}\big(\eta_1-p_1-q_1\big)\bar{x}_2
     =q_2. \end{array}\right.
 \Longleftrightarrow
     \left\{ \begin{array}{ll}
     \eta_1-\left[\bar{x}_2^T\eta_2+\frac{\|d_2\|^2-d_1^2}{\|x_2\|}\right]< 0, \\
     \bar{x}_2^Tp_2-\bar{x}_2^Tq_2+\frac{\|d_2\|^2-d_1^2}{\|x_2\|} = p_1 + q_1, \\
     -q_1\bar{x}_2 = q_2.
     \end{array}\right.\\
& \Longleftrightarrow&
     \left\{ \begin{array}{ll}
     \eta_1-\left[\bar{x}_2^T\eta_2+\frac{\|d_2\|^2-d_1^2}{\|x_2\|}\right]< 0,\\
     \bar{x}_2^Tp_2+\frac{\|d_2\|^2-d_1^2}{\|x_2\|} = p_1, \\
     -q_1\bar{x}_2 = q_2.
     \end{array}\right.
 \Longleftrightarrow
     \left\{ \begin{array}{ll}
     q_1> 0,\\
     \bar{x}_2^Tp_2+\frac{\|d_2\|^2-d_1^2}{\|x_2\|} = p_1, \\
     -q_1\bar{x}_2 = q_2.
     \end{array}\right.\\
 &\Longleftrightarrow&
     \left\{ \begin{array}{ll}
     q\in \Re_{++}\hat{x},\\
     p\in {\rm bd}T_{\cal K}^2(x;d).
     \end{array}\right.
 \end{eqnarray*}
Hence either $p\in T^2_{\K}(x;d)$ and $q=0$ or $ q\in \Re_{++}\hat{x}$ and $p\in {\rm bd}T_{\cal K}^2(x;d)$.
 \noindent

 \medskip
 $ {\bf (iv)}$-3. $d\perp \hat{x}$ and $w\in \Re_{++}\hat{x}$.
 Then
 $ \bar{x}_2^Td_2 = d_1$ and
 $\bar{x}_2^Tw_2= -w_1$.
 Hence
 \begin{equation}\label{eqn09}
  \xi_1-\bar{x}_2^T\xi_2= d_1-w_1-\bar{x}_2^T(d_2-w_2)=-2w_1<0,
  \end{equation}  which implies $\xi\in T_{\K}(x)^c$ by Lemma \ref{Lem3.3}. Thus by Theorem  \ref{projection-second-order diff}(iv),
 $\Pi^{''}_{\K}(z;\xi,\eta)=p$ takes the form
 \begin{equation}\label{add-4-1}
 \left\{\begin{array}{ll} & \frac{1}{2}\left(\eta_1+\bar{x}_2^T\eta_2
 +\frac{\|\xi_2\|^2-(\bar{x}_2^T\xi_2)^2}{\|x_2\|} \right)=p_1,\\
 &\frac{1}{2}\left(\eta_1-\bar{x}_2^T\eta_2-\frac{\|\xi_2\|^2-3(\bar{x}_2^T\xi_2)^2
  +2\xi_1\bar{x}_2^T\xi_2}{\|x_2\|}\right)\bar{x}_2
  +\eta_2+\frac{\xi_1-\bar{x}_2^T\xi_2}{\|x_2\|}\xi_2=p_2.
 \end{array}\right.\end{equation}
 Note that
\begin{eqnarray*}
 \xi_1\bar{x}_2^T\xi_2-(\bar{x}_2^T\xi_2)^2&=& (\xi_1-\bar{x}_2^T\xi_2)\bar{x}_2^T\xi_2=-2w_1(d_1+w_1), \\
  \|\xi_2\|^2-(\bar{x}_2^T\xi_2)^2& =&
  \|d_2-w_2\|^2-(d_1+w_1)^2
 =\|d_2\|^2-d_1^2,
 \end{eqnarray*}
 where we have used the fact $d\perp w$ and $w_2=-w_1\bar{x}_2$ due to $d\perp \hat{x}$ and $w\in \Re_{++}\hat{x}$.
Therefore
 \begin{eqnarray}
\frac{\|\xi_2\|^2-3(\bar{x}_2^T\xi_2)^2
  +2\xi_1\bar{x}_2^T\xi_2}{\|x_2\|}& =&
  \frac{\|\xi_2\|^2-(\bar{x}_2^T\xi_2)^2}{\|x_2\|}+2\frac{\xi_1\bar{x}_2^T\xi_2-(\bar{x}_2^T\xi_2)^2
  }{\|x_2\|}\nonumber\\
  &=&  \frac{\|d_2\|^2-d_1^2}{\|x_2\|}-4\frac{w_1(d_1+w_1)}{\|x_2\|}.\label{eqn4-4-2}
   \end{eqnarray}
   Putting (\ref{eqn09})-(\ref{eqn4-4-2}) into (\ref{add-4-1}) yields
  \begin{eqnarray*}
 &&\Pi_{\K}^{''}(z;\xi,\eta)=p\\
 &\Longleftrightarrow& \left\{\begin{array}{ll} \label{4-2newnnn}
 & \frac{1}{2}\left(\eta_1+\bar{x}_2^T\eta_2+\frac{\|d_2\|^2-d_1^2}{\|x_2\|} \right)=p_1 ,\\
 & \left(\frac{1}{2}\eta_1-\frac{1}{2}\left[\bar{x}_2^T\eta_2+
 \frac{\|d_2\|^2-d_1^2}{\|x_2\|}\right]+2\frac{w_1}{\|x_2\|}(d_1+w_1)\right)\bar{x}_2
 -2\frac{w_1}{\|x_2\|}(d_2-w_2)=q_2.
 \end{array}\right. \\
 &\Longleftrightarrow &
 \left\{\begin{array}{ll} \label{4-2new}
 & \frac{1}{2}\left(\eta_1+\bar{x}_2^T\eta_2+\frac{\|d_2\|^2-d_1^2}{\|x_2\|} \right)=p_1, \\
 & \left(-q_1+2\frac{w_1}{\|x_2\|}(d_1+w_1)\right)\bar{x}_2
 -2\frac{w_1}{\|x_2\|}(d_2-w_2)=q_2.
 \end{array}\right. \\
 &\Longleftrightarrow &
 \left\{\begin{array}{ll} \label{4-2newn}
 & \frac{1}{2}\left(\eta_1+\bar{x}_2^T\eta_2+\frac{\|d_2\|^2-d_1^2}{\|x_2\|} \right)=p_1, \\
 & -q_1\bar{x}_2-2\frac{d_1w_2}{\|x_2\|}-2\frac{w_1d_2}{\|x_2\|}=q_2.
 \end{array}\right.\\
 &\Longleftrightarrow & \left\{\begin{array}{ll} \label{4-2newnn}
 & \bar{x}_2^Tp_2+\frac{\|d_2\|^2-d_1^2}{\|x_2\|}=p_1, \\
 & -q_1\bar{x}_2-2\frac{w_1d_2}{\|x_2\|}-2\frac{d_1w_2}{\|x_2\|}=q_2.
 \end{array}\right.
 \end{eqnarray*}
 where the third equivalence uses the fact $w_2=-w_1\bar{x}_2$ due to $w\in \Re_{++}\hat{x}$ and  the last step follows from substituting  the expression for $q_2$ in the second equation  into the first one to obtain
 \[
 \bar{x}_2^Tq_2=\bar{x}_2^T\left(-q_1\bar{x}_2-2\frac{w_1d_2}{\|x_2\|}-2\frac{d_1w_2}{\|x_2\|} \right)=-q_1.
 \]
  The desired result follows from noting that $p\in {\rm bd}T_{\cal K}^2(x;d)$ if and only if $ \bar{x}_2^Tp_2+\frac{\|d_2\|^2-d_1^2}{\|x_2\|}=p_1$ by virtue of Lemma  \ref{Lem3.3}.

 \medskip
 \noindent
{\bf Case (v)}  $x=0$ and $y\in {\rm bd}\K\backslash \{0\}$.
 The proof is omitted, since  this case is symmetric to Case (iv).

 \medskip
 \noindent
 {\bf Case (vi)} $x=0$ and $y=0$.   Since
 $\Omega$ is cone, according to the definition of second-order tangent
set, we have
 \[
 T^2_{\Omega}\big((0,0);(d,w)\big)=T_{\Omega}(d,w).
 \]
 From all the above, the proof is complete.
 \endproof

 \section{Second-order optimality conditions for SOCMPCC}

In this section, as an application of the second-order tangent set for the SOC complementarity set, we consider second-order  optimality conditions for the mathematical programming with second-order cone complementarity constraints
 (SOCMPCC):
 \begin{equation} \label{SOCMPCC}
 \min   \   f(x)
  \quad {\rm s.t.} \ \ \K \ni G(x) \perp H(x) \in \K,
 \end{equation}
 where $f:\Re^n\to \Re$ and $G,H:\Re^n\to \Re^m$ are second-order continuously differentiable.
 For simplicity, we restrict our attention on the simpler case, i.e.,  $\K$ is a $m$-dimensional  second-order cone. All  analysis
 can be easily carried over to more general cases where $\K$ is a Cartesian product of some second-order cones. SOCMPCC is an important class of optimization problems that has many applications. We refer the reader to \cite{YZ15,YZ16} and the reference within for applications and the first-order necessary optimality conditions.

  Denote by $F(x):=(G(x),H(x))$. Then SOCMPCC (\ref{SOCMPCC}) can be rewritten as
 \begin{equation} \label{SOCMPCCre}
 \min   \   f(x)
  \quad {\rm s.t.} \ \ F(x) \in \Omega.
 \end{equation}
For  a convex set-constrained optimization problem in the form of (\ref{SOCMPCCre}) where $\Omega$ is replaced by  a convex closed set $K$ (see \cite[(3.93)]{BS}), second-order optimality conditions that involve the second-order tangent set to $K$  have been developed in  \cite{BCS99, BS}. In particular when the convex set $K$ is not polyhedral, the second-order tangent set to $K$  is needed in the  second-order optimality conditions.  However, {if set $\Omega$ in problem  (\ref{SOCMPCCre}) is} nonconvex, these optimality conditions are not applicable in general. In what follows, we will establish the second-order optimality conditions
for the SOCMPCC, which  is not a convex set-constrained optimization problem.
We would like to emphasize that, even if the second-order cone complementarity set is nonconvex, its tangent cone and second-order tangent set have nice properties so that some of the theories in the second-order optimality conditions for a convex set-constrained optimization problem still hold.
 This observation relies heavily on the exact formula of tangent cone and second-order tangent set established in the previous section.

First we  present some results needed for further analysis.
Recall that the regular tangent cone is always convex. The following result shows that the regular tangent cone to the SOC complementarity set $\Omega$ is not only convex but is a subspace.
  \begin{proposition}\label{frechet-tangent-cone}
 For any $(x,y)\in \Omega$,
 \begin{eqnarray*}
 \lefteqn{\widehat{T}_\Omega(x,y)={\rm lin} T_{\Omega}(x,y)}\\
 &=&
 \left\{(d,w)\left|
 \begin{array}{lll}
 & d\in \Re^m,\ w=0, \ & {\rm if}\ x\in {\rm int}{\cal K} \ {\rm and}\ y=0; \\
 & d=0,\ w\in \Re^m, \ & {\rm if} \ x=0\ {\rm and}\ y\in {\rm int}{\cal K}; \\
 & x_1\hat{w}-y_1d\in \Re x, \ \ d\perp y,\ w\perp x, \ \
 & {\rm if} \ x, y\in {\rm bd}{\cal K}\backslash\{0\}; \\
 & d\perp \hat{x},\ w=0,  \ \ & {\rm if}\ x\in {\rm bd}{\cal K}\backslash \{0\}\ {\rm and}\ y=0; \\
  & d=0, \ w\perp \hat{y}, \ \ & {\rm if}\ x=0\ {\rm and}\ y\in {\rm bd}{\cal K}\backslash \{0\}; \\
 & d=0,\ w=0, \ \ & {\rm if}\ x=0 \ {\rm and} \ y=0.
 \end{array}\right.\right\}.
 \end{eqnarray*}
 \end{proposition}
 \beginproof
 {The formula of ${\rm lin} T_{\Omega}(x,y)$ is clear from that of $T_{\Omega}(x,y)$ in  Lemma \ref{tangentcone}. According to the tangent-normal polarity as in Lemma \ref{polarity}, we can obtain the formula of $\widehat{T}_\Omega(x,y)$ by taking the polar of the limiting normal cone to $\Omega$ given in \cite[Theorem 5.1]{yzhou}. The obtained formula of ${\rm lin} T_{\Omega}(x,y)$ and $\widehat{T}_\Omega(x,y)$ shows that they have the same expression.} \endproof

 {The exact formula established in Theorem \ref{formula-regular normal cone-1} and Proposition \ref{frechet-tangent-cone} immediately imply the following results.}
 \begin{corollary}\label{cor-lin-tangent}
 For all $(x,y)\in \Omega$,
$
 T_{\Omega}(x,y)+\widehat{T}_{\Omega}(x,y)=T_{\Omega}(x,y).
$
 \end{corollary}



 \begin{proposition}\label{second-order-tangent-cone-add}
 For $(x,y)\in \Omega$ and $(d,w)\in T_{\Omega}(x,y)$,
 \[
 T^2_{\Omega}\big((x,y);(d,w)\big) + \widehat{T}_{\Omega}(x,y)
 = T^2_{\Omega}\big((x,y);(d,w)\big).
 \]
 \end{proposition}
 \beginproof
  {The inclusion ``$\supseteq$" is clear, since $(0,0)\in \widehat{T}_{\Omega}(x,y)$.
For all cases except where  $x, y \in {\rm bd}\K\backslash\{0\}$,  it is easy to see that  ``$\subseteq$" can be achieved by using the formula of $T^2_{\Omega}\big((x,y);(d,w)\big)$ given in Theorem \ref{formula-regular normal cone-1} and the formula of $\widehat{T}_\Omega(x,y)$ given in Proposition \ref{frechet-tangent-cone}. Now consider the case where  $x, y \in {\rm bd}\K\backslash\{0\}$.
Let $(p,q)\in T^2_{\Omega}\big((x,y);(d,w)\big)$ and $(u,v)\in \widehat{T}_{\Omega}(x,y)$.  Since $p\in {\rm bd}T_{\K}^2(x;d)$, then $\hat{x}^Tp= \|d_2\|^2-d_1^2$
 by Lemma \ref{Lem3.3}. Hence $\hat{x}^T(p+u)=\hat{x}^Tp+\hat{x}^Tu= \|d_2\|^2-d_1^2$ due to the fact $u\perp \hat{x}$ (since $u\perp y$ by Proposition \ref{frechet-tangent-cone} and $y\in \Re_{++}\hat{x}$ ). This means $p+u\in {\rm bd}T_{\K}^2(x;d)$. Similarly, we can obtain $q+v\in {\rm bd}T_{\K}^2(y;w)$. Since $(u,v)\in \widehat{T}_{\Omega}(x,y)$, it follows from Proposition \ref{frechet-tangent-cone} that there exists $\tau\in \Re$ such that $x_1\hat{v}-y_1u=\tau x$. Thus
 \begin{equation}\label{re-eq-a11a}
 x_1v_2+y_1u_2=-\tau x_2=-\frac{x_1v_1-y_1u_1}{x_1}x_2=
 -v_1x_2-u_1y_2,
 \end{equation}
 where the last step comes from Lemma \ref{lem2.1}. Since $x, y \in {\rm bd}\K\backslash\{0\}$,  $(p,q)\in T^2_{\Omega}\big((x,y);(d,w)\big)$, it follows from Theorem \ref{formula-regular normal cone-1} that
 \begin{equation}
 (p,q)\in T^2_{\Omega}\big((x,y);(d,w)\big) \Longleftrightarrow  \left\{  \begin{array}{l}
 p\in {\rm bd}T_{\K}^2(x;d), q \in {\rm bd}T_{\K}^2(y;w) \\
 \xi-p_1y_2-q_1x_2=x_1q_2+y_1p_2\end{array} \right. , \label{xi}
 \end{equation}
 where $ \xi:=  (x_1w_1-y_1d_1)\left(\frac{{w_2-w_1\bar{y}_2}}{y_1}-\frac{d_2-d_1\bar{x}_2}{x_1}\right).
$ This, together with (\ref{re-eq-a11a}), implies
  \begin{eqnarray*}
 x_1(q_2+v_2)+y_1(p_2+u_2)&=& x_1q_2+y_1p_2-v_1x_2-u_1y_2\\
 &=& \xi-p_1y_2-q_1x_2-v_1x_2-u_1y_2\\
 &=& \xi-(p_1+u_1)y_2-(q_1+v_1)x_2.
 \end{eqnarray*}
 Hence together with  $p+u\in {\rm bd}T_{\K}^2(x;d)$ and $q+v\in {\rm bd}T_{\K}^2(y;w)$, we have that  $(p+u,q+v)\in T^2_{\Omega}\big((x,y);(d,w)\big)$ by virtue of (\ref{xi}).
 \endproof

 With these preparations, we are now ready to develop a second-order necessary optimality condition for SOCMPCCs. Define the Lagrange function as
$L(x,\lambda):= f(x)+\langle F(x), \lambda\rangle$ and the following three multiplier  sets
\begin{eqnarray*}
 \Lambda^c(x^*) &:=& \{\lambda \, | \, \nabla_x L(x^*,\lambda)=0,
      \ \lambda\in N^c_{\Omega}(F(x^*))\} , \\
  \Lambda (x^*) &:=& \{\lambda \, | \, \nabla_x L(x^*,\lambda)=0,
      \ \lambda\in {N}_{\Omega}(F(x^*))\},\\
      \Lambda^F(x^*) &:=& \{\lambda \, | \, \nabla_x L(x^*,\lambda)=0,
      \ \lambda\in \widehat{N}_{\Omega}(F(x^*))\}.
 \end{eqnarray*}
Denote by $
 C(x^*):=\big\{d|\, \nabla f(x^*)d\leq 0, \nabla F(x^*)d\in T_{\Omega}(F(x^*)) \big\}
$  the critical cone.
Note that if there exists $\lambda\in \Lambda^F({x^*})$, then $C(x^*)=\big\{d|\, \nabla f(x^*)d= 0, \nabla F(x^*)d\in T_{\Omega}(F(x^*)) \big\}$.

 \begin{theorem}\label{dual-conditions}
Let $x^*$ be a locally optimal solution of SOCMPCC.
Suppose that the non-degeneracy condition
 \begin{equation}\label{conditions-revision}
 {\nabla} F(x^*)\Re^n+{\rm lin}T_{\Omega}(F(x^*)) = \Re^{2m}
 \end{equation} holds. Then
  $\Lambda^c(x^*) =\Lambda(x^*)=\Lambda^F(x^*)=\{\lambda_0\}$ and
 \[
 \nabla^2_{xx}L(x^*,\lambda_0)(d,d) - \sigma \left( \lambda_0| T^2_{\Omega}
 (F(x^*);\nabla F(x^*)d) \right) \geq 0, \quad
 \forall d\in C(x^*).
 \]
 \end{theorem}
 \beginproof
%
 \noindent
 Step 1. We prove $\Lambda^c(x^*) =\Lambda(x^*)=\Lambda^F(x^*)=\{\lambda_0\}$.
 Since
 $
 \nabla F(x^*)\Re^n+{\rm lin}T_{\Omega}(F(x^*))=
 \nabla F(x^*)\Re^n+\widehat{T}_{\Omega}(F(x^*))=\Re^{2m},$ taking polars on the both sides of the above equation, by the rule for polar cones \cite[Corollary 11.25]{RW98} and the fact that
 $(\widehat{T}_{\Omega})^{\circ}=N^c_{\Omega}$,
 we have
 \begin{equation} \label{clarke-normal}
 \nabla F(x^*)^T\lambda=0, \  \lambda\in N^c_{\Omega}(F(x^*))
 \quad \Longrightarrow \quad \lambda = 0.
 \end{equation}
 Suppose that  $\lambda^1, \lambda^2\in \Lambda^c(x^*)$.
 Then $\lambda^1-\lambda^2$ satisfies $\nabla F(x^*)^T(\lambda^1-\lambda^2)=0$
 and $\lambda^1-\lambda^2\in N^c_{\Omega}(F(x^*))$ since  $N^c_{\Omega}(F(x))$
 is subspace (because $\widehat{T}_{\Omega}(x,y)$ is a subspace and
 $N^c_{\Omega}=(\widehat{T}_{\Omega})^{\circ}$). Thus
 $\lambda^1=\lambda^2$ by (\ref{clarke-normal}). This means that
 $\Lambda^c(x^*)$ is a singleton. Since
 $\Lambda^F(x^*) \subseteq \Lambda(x^*) \subseteq \Lambda^c(x^*)$, it remains to
 show that $\Lambda^F(x^*)$ is nonempty. Since
 $N_{\Omega}\subseteq N^c_{\Omega}$, the condition
 (\ref{clarke-normal}) ensures
 \begin{equation}
 \nabla F(x^*)^T\lambda=0, \  \lambda\in N_{\Omega}(F(x^*))
 \quad \Longrightarrow \quad \lambda=0,\label{MR}
 \end{equation}
  which in turn implies that the system $F(x)-\Omega$ is metrically regular
 at $(x^*,0)$. Thus according to \cite[Theorem 4]{GO16}, Proposition \ref{frechet-tangent-cone}, and Corollary \ref{cor-lin-tangent}, we have
 $\widehat{N}_S(x)=\nabla F(x)^T\widehat{N}_{\Omega}(F(x))$, where
 $S:=\{x \, | \, F(x)\in \Omega\}$. As $x^*$ is a local optimal
 solution of problem (\ref{SOCMPCCre}), we have
 \(
 0\in \nabla f(x^*)^T+\widehat{N}_{S}(x^*)
 =\nabla f(x^*)^T+\nabla F(x^*)^T\widehat{N}_{\Omega}(F(x^*))
 \),
 which indicates that $\Lambda^F(x^*)$ is nonempty. Hence $\Lambda^c(x^*), \Lambda(x^*),\Lambda^F(x^*)$
are all singleton and coincide with each other. Let us denote the unique element by $\lambda_0$.

 \medskip
 \noindent
 Step 2. We show that for all $d\in C(x^*)$ and for any convex subset
 ${\cal T}(d)$ in $T^2_{\Omega}(F(x^*);\nabla F(x^*)d)$,
$
 \nabla^2_{xx}L(x^*,\lambda_0)(d,d)-\sigma \left( \lambda_0|{\cal T}(d) \right)\geq 0.
 $
  {The idea of the proof is inspired by the arguments in \cite[Theorem 3.1]{BCS99} and using the properties of tangent cone and second-order tangent set discussed above.} {For the sake of completeness, we give the detailed proof  here.}
  Consider the set
 $\Gamma(d):={\rm cl}\{{\cal T}(d)+\widehat{T}_{\Omega}(F(x^*))\}$. {Since the regular tangent cone is convex, the set $\Gamma(d)$}  is
 closed  and convex. Moreover, it follows from Proposition \ref{second-order-tangent-cone-add}
 and the fact that the second-order tangent set is closed that  $\Gamma(d)\subseteq T^2_{\Omega}(F(x^*);\nabla F(x^*)d)$.
 Because $x^*$ is locally optimal of problem (\ref{SOCMPCCre}),  {by definition of the second-order tangent cone,  we can show that
 $$\nabla f(x^*)w +\nabla^2f(x^*)(d,d)  \geq 0, \quad \forall d\in C(x^*), w\in T^2_C\big(x^*;d
 \big),$$
 where $C$ denotes the feasible region of problem (\ref{SOCMPCCre}).  Since (\ref{MR}) holds, by  \cite[Proposition 13.13]{RW98},  the chain rule for tangent sets (\ref{chainruletangentset}) holds with $\Theta$ taken as $\Omega$.}  It follows that for all $d\in C(x^*)$,
 the following optimization problem
 \[
 \left\{
 \begin{array}{ll}
 \displaystyle \min_w  & \nabla f(x^*)w+\nabla^2f(x^*)(d,d)  \\
 {\rm s.t.} & \nabla F(x^*)w + \nabla^2F(x^*)(d,d)
              \in T^2_\Omega\big(F(x^*);\nabla F(x^*)d \big)
 \end{array}
 \right.
 \]
 has nonnegative optimal value.
 Since $\Gamma(d)\subseteq T^2_{\Omega}(F(x^*);\nabla F(x^*)d)$, it
 is clear that the following convex set constrained problem
 \begin{equation} \label{convex-subproblem}
 \left\{
 \begin{array}{ll}
\displaystyle \min_w  & \nabla f(x^*)w+\nabla^2f(x^*)(d,d)  \\
 {\rm s.t.} & \nabla F(x^*)w + \nabla^2F(x^*)(d,d)
              \in \Gamma(d)
 \end{array}
 \right.
 \end{equation}
 has nonnegative optimal value as well. Since the optimization problem  (\ref{convex-subproblem}) can be put into the form of problem \cite[(2.291)]{BS} involving an indicator function of set $\Gamma(d)$ and  the dual problem of \cite[(2.291)]{BS} is in the form of \cite[(2.298)]{BS} and  the conjugate function of an indicator function is the support function, the dual problem of  (\ref{convex-subproblem}) is
 $$ \max_{\lambda} \left \{ \inf_{w} {\cal L}(w, \lambda)
 -\sigma \left( \lambda|\Gamma(d) \right) \right \},$$
 where ${\cal L}(w,\lambda):=\nabla_x L(x^*,\lambda)w +\nabla_{xx}^2 L(x^*, \lambda)(d,d) $ is the Lagrange function of (\ref{convex-subproblem}).
  {Note that
 \[\sigma(\lambda|\Gamma(d))=\sigma\big(\lambda|{\cal T}(d)+\widehat{T}_{\Omega}(F(x^*))\big)=\sigma(\lambda|{\cal T}(d))+
 \sigma\big(\lambda|\widehat{T}_{\Omega}(F(x^*))\big)=+\infty,\]
 whenever
 $\lambda\notin [\widehat{T}_{\Omega}(F(x))]^\circ=N^c_{\Omega}(F(x))$.}
  Therefore, the dual problem of (\ref{convex-subproblem}) is
 \begin{equation}\label{dual-problem}
 \max_{\lambda\in \Lambda^c(x^*)} \left \{ \nabla^2_{xx}L(x^*,\lambda)(d,d)
 -\sigma \left( \lambda|\Gamma(d) \right) \right \}
 = \nabla^2_{xx}L(x^*,\lambda_0)(d,d)-\sigma \left( \lambda_0|\Gamma(d) \right),
 \end{equation}
 where the equality holds since $\Lambda^c(x^*)=\{\lambda_0\}$ by Step 1.

  Since ${\rm lin}T_{\Omega}(F(x^*))=\widehat{T}_{\Omega}(F(x^*))$ by Proposition \ref{frechet-tangent-cone} and
 ${\rm lin}T_{\Omega}(F(x^*))$ is a subspace, we have ${\rm lin}T_{\Omega}(F(x^*))=-\widehat{T}_{\Omega}(F(x^*)).$
 Hence  condition (\ref{conditions-revision}) is
 $\nabla F(x^*)\Re^n-\widehat{T}_{\Omega}(F(x^*))=\Re^{2m}$, which in turn implies
 $\nabla F(x^*)\Re^n-\big({\cal T}(d)+\widehat{T}_{\Omega}(F(x^*))\big)=\Re^{2m}$. Hence
 $\nabla F(\bar x)\Re^n-\Gamma(d)=\Re^{2m}$. So
the
 Robinson's  constraint qualification
 (see \cite[(2.313)]{BS}) for problem (\ref{convex-subproblem}) holds. It ensures that the zero dual gap property holds (see
\cite[Theorem 2.165]{BS}). Hence the optimal value of the dual
 problem (\ref{dual-problem}) is equal to the optimal value of  problem  (\ref{convex-subproblem}) and hence nonnegative. In addition, noting
 ${\cal T}(d)\subseteq \Gamma(d)$,
 $\sigma(\lambda_0| {\cal T}(d))\leq \sigma(\lambda_0|\Gamma(d))$, which further implies that
 \begin{equation} \label{second-order-conditions}
  \nabla^2_{xx}L(x^*,\lambda_0)(d,d)-\sigma \left( \lambda_0|{\cal T}(d) \right) \geq 0.
 \end{equation}

 \noindent
 Step 3. Note that $T^2_{\Omega}(F(x^*);\nabla F(x^*)d)=\bigcup_{a\in T^2_{\Omega}(F(x^*);\nabla F(x^*)d)}\{a\}$ is the union of
 convex sets. For each $a\in T^2_{\Omega}(F(x^*);\nabla F(x^*)d)$, by (\ref{second-order-conditions}) we have
 $$  \nabla^2_{xx}L(x^*,\lambda_0)(d,d)-\langle \lambda_0,a \rangle  \geq 0.$$
It then yields the desired result
 \[
 \nabla^2_{xx}L(x^*,\lambda_0)(d,d) - \sigma \left( \lambda_0|T^2_{\Omega}(F(x^*);
 \nabla F(x^*)d) \right) \geq 0.
 \]
 \endproof

 \begin{remark}
 The nondegeneracy condition (\ref{conditions-revision}), together with the special geometric structure of second-order cone complementarity set, can ensure not only the uniqueness of Lagrangian multiplier in Step 1, but also the zero-dual gap property between (\ref{convex-subproblem}) and (\ref{dual-problem}) in Step 2. The nondegeneracy condition, stronger than the Robinson's constraint qualification, is a generalization of linear independence constraint qualification in the conic case. We refer to \cite[Proposition 4.75]{BS} for the detailed discussion on the relationship between nondegeneracy condition and uniqueness of multiplier in the convex case.
 \end{remark}


 We next derive the exact formula for the support function of the second-order tangent set to the SOC complementarity set needed in applying Theorem \ref{dual-conditions}. Under the assumption of Theorem
\ref{dual-conditions} we have $C(x^*)=\big\{d|\, \nabla f(x^*)d= 0, \nabla F(x^*)d\in T_{\Omega}(F(x^*)) \big\}$. Thus $d\in C(x^*)$ if and only if $\nabla F(x^*)d\in T_{\Omega}(F(x^*))$ and $\langle \lambda_0, \nabla F(x^*)d\rangle =0$. Therefore  the following results will be useful.

 \begin{proposition}\label{formula-regular normal cone}
 { For $(x,y)\in \Omega$ and $(d,w)\in T_\Omega(x,y)$, take $(u,v)\in \widehat{N}_{\Omega}(x,y)$ such that $\langle (u,v),(d,w) \rangle=0$. Then }
 \[
 \sigma\big((u,v)|T^2_{\Omega}((x,y);(d,w))\big)=\left\{\begin{array}{lll}
 0, \ & {\rm if} \ x\in {\rm int} \K\ {\rm and}\ y=0; \\
 0, \ & {\rm if}\ x=0 \ {\rm and}\ y \in {\rm int}{\cal K};\\
 0, \  & {\rm if}\ x=0 \ {\rm and} \ y=0.
 \end{array}\right\}.
 \]
 If $x\in {\rm bd}\K\backslash \{0\}$ and $y=0$, then
 \begin{eqnarray*}
 &&\sigma\big((u,v)|T^2_{\Omega}((x,y);(d,w))\big)\\
 &=&\left\{\begin{array}{lll}
 0,  \ & {\rm if} \ d\in  {\rm int}T_{\K}(x)\ {\rm and}\ w=0;
 \\
 -\frac{u_1}{x_1}(d_1^2-\|d_2\|^2)-2\frac{w_1d_2^Tv_2}{\|x_2\|}-2\frac{d_1w_2^Tv_2}{\|x_2\|}, \ \ & {\rm if}\ d\in {\rm bd}T_{\K}(x) \ {\rm and} \ w\in \Re_{+}\hat{x}.
 \end{array}\right\}.
 \end{eqnarray*}
 If $x=0$ and $y\in {\rm bd}\K\backslash \{0\}$, then
 \begin{eqnarray*}
 &&\sigma\big((u,v)|T^2_{\Omega}((x,y);(d,w))\big)\\
 &=&\left\{\begin{array}{lll}
 0,  \ & {\rm if} \ d=0\ {\rm and}\ w\in  {\rm int}T_{\K}(y); \\
 -\frac{v_1}{y_1}(w_1^2-\|w_2\|^2)-2\frac{d_1w_2^Tu_2}{\|y_2\|}
 -2\frac{w_1d_2^Tu_2}{\|y_2\|}, \ \ & {\rm if}\  d\in \Re_{+}\hat{y} \ {\rm and} \  w\in {\rm bd}T_{\K}(y).
  \end{array}\right\}.
 \end{eqnarray*}
 If $x,y\in {\rm bd}\K\backslash \{0\}$, then
 \begin{eqnarray*}
 &&\sigma\big((u,v)|T^2_{\Omega}((x,y);(d,w))\big)\\
 &=& \frac{x_1u_1+y_1v_1}{x^2_1}\Big(\|d_2\|^2-d_1^2\Big)+\frac{x_1w_1-y_1d_1}{x_1y_1} \Big(w^Tv-d^Tu\Big).
  \end{eqnarray*}
 \end{proposition}

 \beginproof {For $(x,y) \in \Omega$, take $(d,w)\in T_\Omega(x,y)$,  $(p,q)\in T^2_\Omega((x,y);(d,w))$} with the exact formula given in Theorem \ref{formula-regular normal cone-1} and  $(u,v) \in \widehat{N}_\Omega(x,y) $ whose exact formula can be found  in \cite[Theorem 3.1]{yzhou}.

 \noindent {\bf Case (i)} $x\in {\rm int}\K$ and $y=0$.  In this case $u=0$ and $q=0$. Hence
 $\sigma\big((u,v)|T^2_{\Omega}((x,y);(d,w))\big)=
 \max\{\langle (u,v),(p,q)\rangle|(p,q)\in T_\Omega^2((x,y);(d,w))\}=0$. The proof for the  case of $x=0$ and $y\in {\rm int}\K$ is similar and hence we omit it.

  \noindent {\bf Case (ii)}  $x\in {\rm bd}\K\backslash \{0\}$ and $y=0$. Then $u\in \Re_-\hat{x}$ and $v\in \hat{x}^\circ$ by the formula of $\widehat{N}_\Omega(x,y)$.

  ${\bf  (ii)}$-1.  Suppose further that $d\in {\rm int} T_{\K}(x)\ {\rm and}\ w=0$. Then by Theorem \ref{formula-regular normal cone-1}, $q=0$. Since $0=\langle u,d \rangle+ \langle v,w \rangle =\langle u,d \rangle $, which together with the fact that $d\in {\rm int} T_{\K}(x)$ (i.e., $d^T\hat{x}>0$) implies $u=0$.
  Hence $\sigma\big((u,v)|T^2_{\Omega}((x,y);(d,w))\big)=\max\{\langle (u,v), (p,q)\rangle|(p,q)\in T^2_{\Omega}((x,y);(d,w))\}=0.$

  ${\bf  (ii)}$-2.  Suppose further that $d\in {\rm bd}T_{\K}(x)\ {\rm and } $ $ w=0$,  then $q=0$ or $q\in \Re_+\hat{x}$ by Theorem \ref{formula-regular normal cone-1}. Hence
   \begin{eqnarray*}
   \sigma\big((u,v)|T^2_{\Omega}\big((x,y);(d,w)\big)\big)
   &=&\max\Big\{\sigma(u| T^2_{\K}(x;d)),\sigma(u|{\rm bd}T^2_\K(x;d))+\sigma(v|\Re_+\hat{x})\Big\} \\
   &=& \sigma(u| T^2_{\K}(x;d))=  -\frac{u_1}{x_1}(d_1^2-\|d_2\|^2),
   \end{eqnarray*}
   where the second equality holds because
   $\sigma(v|\Re_+\hat{x})=0$ since $v\in \hat{x}^\circ$, and the last step comes from the fact that  since $u\in \Re_-\hat{x}$,
    $\langle u,p \rangle=\frac{u_1}{x_1} \langle \hat{x}, p \rangle\leq \frac{u_1}{x_1}(\|d_2\|^2-d_1^2)$ for all $p\in T^2_{\K}(x;d)$ by Lemma \ref{Lem3.3}, and the maximum can be attained by letting $p=\frac{\|d_2\|^2-d_1^2}{2x_1^2}\hat{x}$.

{Now consider the case where $d\in {\rm bd}T_{\K}(x)\ {\rm and } $ and $w\in \Re_{++}\hat{x}$}. From the formula for $  {\rm bd}T_{\K}(x)$ in this case, we get  $d\perp \hat{x}$.  Hence  $\langle v, w\rangle =\langle (u,v),(d,w) \rangle=0$ taking into the account that $u\in \Re_-\hat{x}$. It further implies that $v\perp \hat{x}$ (i.e., $v_1=\bar{x}_2^Tv_2$), because $w\in \Re_{++}\hat{x}$.
   Hence
 \begin{eqnarray*}\sigma\big((u,v)|T^2_{\Omega}\big((x,y);(d,w)\big)\big)
 &=&\sigma(u| {\rm bd}T^2_{\K}(x;d)) +\langle v,q\rangle \\
 &=&-\frac{u_1}{x_1}(d_1^2-\|d_2\|^2)+v_1q_1-q_1v_2^T\bar{x}_2-2\frac{w_1d_2^Tv_2}{\|x_2\|}-2\frac{d_1w_2^Tv_2}{\|x_2\|}\\
 &=& -\frac{u_1}{x_1}(d_1^2-\|d_2\|^2)-2\frac{w_1d_2^Tv_2}{\|x_2\|}-2\frac{d_1w_2^Tv_2}{\|x_2\|}.
 \end{eqnarray*}

 \noindent {\bf Case (iii)} $x=0$ and $y\in {\rm bd}\K\backslash \{0\}$. The argument is similar to the above case.

 \noindent {\bf Case (iv)} $x,y\in {\rm bd}\K\backslash \{0\}$.  Note that in this case since $x^Ty=0$, we have $y =k \hat{x}$ with $k:=y_1/x_1$. Since $(u,v)\in \widehat{N}_{\Omega}(x,y)$ and $(d,w)\in T_{\Omega}(x,y)$,  by the formula of $\widehat{N}_{\Omega}(x,y)$ and $T_{\Omega}(x,y)$ we have  $v\perp y$, $d\perp y$, and there exist $\beta,\gamma\in \Re$ such that
 $
 \hat{u}+kv=\beta x$ and $\hat{w}-kd=rx
 $.  To simplify the notation, let
 \[\xi:=  (x_1w_1-y_1d_1)\left(\frac{w_2+w_1\bar x_2}{y_1}-\frac{d_2-d_1\bar x_2}{x_1}\right).
 \]
 Since $v\perp y$, $y \in \Re \hat{x}$  and $x_1=\|x_2\|\not =0$, we have $v_1-{\bar x}_2^Tv_2=0$. It follows that
  \begin{equation}\label{EE}
 v_2^T\xi = (x_1w_1-y_1d_1)\left(\frac{w_2^Tv_2+w_1v_1}{y_1}-\frac{d_2^Tv_2-d_1v_1}{x_1}\right)
 = \frac{x_1w_1-y_1d_1}{y_1}(w^Tv-d^Tu),
 \end{equation}
 where in the last step we used the fact that
 $d^T\hat{v}=(1/k)d^T(\beta\hat{x}-u)=-(1/k)d^Tu$ since $d\perp \hat{x}$. By the formula of $T^2_\Omega((x,y);(d,w))$ in Theorem \ref{formula-regular normal cone-1} for this case we have
\begin{equation}\label{SOT}
p\in {\rm bd} T^2_{\K}(x;d),\ \ q\in {\rm bd}T^2_{\K}(y;w),\ \
 \xi-
 p_1y_2-q_1x_2=x_1q_2+y_1p_2.\end{equation}
Therefore,
 \begin{eqnarray*}
 &&\langle u,p \rangle + \langle v,q \rangle \\
 &&=\langle \hat{u},\hat{p} \rangle + \langle v,q \rangle =\langle \beta x-kv, \hat{p} \rangle + \langle v,q \rangle = \beta \langle \hat{x}, p \rangle + \langle v,q-k\hat{p} \rangle\\
 &&=\beta\langle \hat{x},p \rangle +v_1 (q_1-kp_1)+v_2^T
 \left( \frac{\xi}{x_1}+(kp_1-q_1)\bar x_2\right) \\
 &&= \beta\langle \hat{x},p \rangle +\frac{1}{x_1}v_2^T\xi
 = \frac{x_1u_1+y_1v_1}{x^2_1}\big(\|d_2\|^2-d_1^2\big)+\frac{x_1w_1-y_1d_1}{x_1y_1} (w^Tv-d^Tu),
 \end{eqnarray*}
 where the forth equality holds by virtue of (\ref{SOT}), the fifty equality holds because $v_1=v_2^T\bar{x}_2$ and the sixth equality holds due to (\ref{EE}) and (\ref{SOT}). The desired formula follows.

 \noindent {\bf Case (v)} $x=0$ and $y=0$. In this case  $ d, w \in {\cal K}, d\perp w $, $(u,v)\in \widehat{N}_{\Omega}(x,y)=(-\K,-\K)$ and $(p,q) \in T^2_{\Omega}\big((x,y);(d,w)\big)=T_\Omega (d,w)$.

 ${\bf (v)}$-1. $d=0$ and $w\in {\rm int}\K$. Since $\langle v,w \rangle =\langle (u,v), (d,w) \rangle=0$ and $v\in -\K$, we have $v=0$. Since $d=0$ and $w\in {\rm int}\K$,  $(p,q)\in T^2_{\Omega}\big((x,y);(d,w)\big)=T_{\Omega}(d,w)$ implies that
 $p=0$. Hence $\langle (u,v), (p,q) \rangle=0$. It follows that $\sigma\big((u,v)|T^2_{\Omega}((x,y);(d,w))\big)=0$.

 ${\bf (v)}$-2. $d\in {\rm int}\K$ and $w=0$. It is similar to the above case.

 ${\bf (v)}$-3. $d,w\in {\rm bd}\K\backslash\{0\}$. Then since $\langle (u,v), (d,w) \rangle=0$ and $(u,v)\in (-\K,-\K)$, we have $u\in \Re_-\hat{d}=\Re_-w$ and $v\in \Re_-\hat{w}=\Re_-d$.
 Since  $(p,q)\in T_{\Omega}(d,w)$ and $d,w\in {\rm bd}\K\backslash\{0\}$, we have $p\perp w$ and $q\perp d$. Hence $p\perp u$ and $q\perp v$. So
 $\langle (u,v), (p,q) \rangle=0$. It follows that $\sigma\big((u,v)|T^2_{\Omega}((x,y);(d,w))\big)=0$.

 ${\bf (v)}$-4. $d=0$ and $w\in {\rm bd}\K\backslash\{0\}$.  Since $\langle v,w \rangle =\langle (u,v), (d,w) \rangle=0$ and $v\in -\K$, we have  $v\in \Re_-\hat{w}$.
 In this case since $(p,q)\in T_{\Omega}(d,w)$ with $d=0$ and $w\in {\rm bd}\K\backslash\{0\}$, we have either   $p=0$ and $q\in T_{\K}(w)$ or $p\in \Re_+\hat{w}$ and $q\perp \hat{w}$. If $p=0$ and $q\in T_{\K}(w)$ (i.e., $\hat{w}^Tq\geq 0$), then
 $\langle (u,v), (p,q) \rangle=\langle v,q \rangle \leq 0$ and the maximum is $0$ which can be attained by letting $q=0$. If $p\in \Re_+\hat{w}$ and $q\perp \hat{w}$, then
 $\langle (u,v), (p,q) \rangle=\langle u,p \rangle \leq 0$, where the last step is  due to $u\in -\K$ and $p\in \Re_+\hat{w}\in \K$, and the maximum is $0$ which can be attained by letting $p=0$. It follows that $\sigma\big((u,v)|T^2_{\Omega}((x,y);(d,w))\big)=0$.

 ${\bf (v)}$-5. $d\in {\rm bd}\K/\{0\}$ and $w=0$. It is similar to the above case by symmetric argument.

 ${\bf (v)}$-6. $d=0$ and $w=0$. In this case  $(p,q)\in T_{\Omega}(d,w)=\Omega$. Since $(u,v)\in (-\K,-\K)$, we have
  $\langle (u,v), (p,q) \rangle \leq 0$ and the maximum is $0$ which can be attained by letting $(p,q)=(0,0)$. It follows that $\sigma\big((u,v)|T^2_{\Omega}((x,y);(d,w))\big)=0$.
 \endproof

 \begin{example}
{Consider the following SOCMPCC. }
\[\begin{array}{ll}
 \min       & f(x):=-x_2^2+x_1-x_4  \\
 {\rm s.t.} & \K \ni G(x):=(x_1,x_3-x_1,x_1-x_2) \perp (-x_2+1,x_1,x_1-x_4)=:H(x) \in \K.
 \end{array}
 \]
\end{example}
  Since $x_1\geq x_1-x_2$ and $-x_2+1\geq 0$, we have $x_2\in [0,1]$. Hence $x_2^2\leq x_2$. Since
 $-x_2+1\geq -x_1+x_4$,  $-x_2+x_1-x_4\geq -1$.
 Thus
$
 -x_2^2+x_1-x_4\geq -x_2+x_1-x_4\geq -1.
$
Hence $x^*=(0,0,0,1)$ is an optimal solution, and $G(x^*)=(0,0,0)$ and $H(x^*)=(1,0,-1)\in {\rm bd}\K\backslash \{0\}$.

Note that
 \[
 \nabla G(x^*)=\begin{bmatrix}
 1 & 0 & 0 & 0\\
  -1 & 0 & 1 & 0 \\
  1 & -1 & 0 & 0
 \end{bmatrix},\ \ \ \ \ \nabla H(x^*)=\begin{bmatrix}
 0 & -1 & 0 & 0\\
 1 & 0 & 0 & 0\\
 1 & 0 & 0 & -1
 \end{bmatrix}\]
 and by the formula of the tangent cone in Lemma \ref{tangentcone}, we have
\begin{equation}\label{Tconeeg} T_{\Omega}(G(x^*),H(x^*))=\left \{ (d,w)\left| \begin{array}{l}
 \mbox{either } d=0, -w_1-w_3\leq 0 \\
 \mbox {or } d=t(1,0,1) \mbox{ for some } t\geq 0, w_1+w_3=0
 \end{array} \right. \right\}.\end{equation}
 It follows that
 \[ {\rm lin}T_{\Omega}(G(x^*),H(x^*))=\{((0,0,0),(\tau_1,\tau_2,-\tau_1))|\tau_1,\tau_2\in \Re\}.
 \]
  For any $v\in \Re^6$, take $\xi=(v_1,v_1-v_3,v_1+v_2,v_3-v_4-v_6)\in \Re^4$ and $\tau=(v_1-v_3+v_4,v_5-v_1,-(v_1-v_3+v_4))\in \Re^3$, then \[
v= \left(\begin{matrix}
 \nabla G(x^*)\\
 \nabla H(x^*)
 \end{matrix}\right)\xi+\left(\begin{matrix}
 0\\
 \tau
 \end{matrix}\right)\in \nabla F(x^*)\Re^4+{\rm lin}T_{\Omega}(F(x^*)).
 \]
 Since $v$ is arbitrarily taken from $\Re^6$,   condition (\ref{conditions-revision}) holds.

 The Lagrangian multiplier system is
 \[\left\{\begin{array}{cc}
 \left(\begin{matrix}
   1 \\ 0 \\ 0 \\-1
 \end{matrix}\right)+
 \lambda^G_1\left(\begin{matrix}
   1 \\ 0 \\ 0 \\0
 \end{matrix}\right)+ \lambda^G_2\left(\begin{matrix}
   -1 \\ 0 \\ 1 \\0
 \end{matrix}\right) + \lambda^G_3
 \left(\begin{matrix}
   1 \\ -1 \\ 0 \\0
 \end{matrix}\right)
 + \lambda^H_1\left(\begin{matrix}
   0 \\ -1 \\ 0 \\0 \end{matrix}\right)+ \lambda^H_2\left(\begin{matrix}
   1 \\ 0 \\ 0 \\0 \end{matrix}\right)+\lambda^H_3\left(\begin{matrix}
   1 \\ 0 \\ 0\\ -1 \end{matrix}\right)=\left(\begin{matrix}
   0 \\ 0 \\ 0 \\0
 \end{matrix}\right), \\
    (\lambda^G,\lambda^H)\in N_{\Omega}\big(G(x^*),H(x^*)\big).
    \end{array}
    \right.
    \]
 Since $G(x^*)=(0,0,0)$ and $H(x^*)=(1,0,-1)\in {\rm bd}\K\backslash \{0\}$, we obtain the  following expression of the limiting normal cone from \cite[Theorem 5.1]{yzhou}
 \begin{flalign*}
 &N_{\Omega}\big(G(x^*),H(x^*)\big)\\
 &=\{(u,v)|v=0 \ \ {\rm or} \ \ u_1+u_3=0, v=t(1,0,1),t\in \Re\ \  {\rm or} \ \ u_1+u_3\leq 0, v=t(1,0,1), t\leq 0\}.
 \end{flalign*}
 Hence the only multipliers $(\lambda^G,\lambda^H)$ satisfying the Lagrangian multiplier system is  $\lambda^G=(-1,0,1)$ and $\lambda^H=(-1,0,-1)$. Note that
 \begin{eqnarray*}
 C(x^*)&=&\{d\in \Re^4| (d_1,d_3-d_1,d_1-d_2,-d_2,d_1,d_1-d_4)\in T_{\Omega}(G(x^*),H(x^*)), d_1\leq d_4\} \\
 &=& \{d=(t,0,t,t)|t\geq 0\},
 \end{eqnarray*}
 where the second equality follows from (\ref{Tconeeg}).

Since $\nabla G(x^*)d=(t,0,t), \nabla H(x^*)d =(0,t,0)$ for any $d=(t,0,t,t)$ with $t\geq 0$ in $C(x^*)$,
  by
   Proposition \ref{formula-regular normal cone} we obtain
   \[
  \sigma\big((\lambda^G,\lambda^H)|T^2_\Omega(G(x^*),H(x^*);\nabla G(x^*)d,\nabla H(x^*)d)\big)=
  -t^2=-d_1^2.
  \]
Since  $
 \nabla^2_{xx}  L(x, \lambda)=\nabla^2f(x)=\begin{bmatrix}
 0 & 0 & 0 & 0\\
  0 & -2 & 0 & 0 \\
  0 & 0 & 0 & 0\\
  0&0&0&0
 \end{bmatrix}$, we have
    \begin{equation}\label{re-eq-2}
    \nabla^2_{xx}L(x^*,\lambda)(d,d)=0,  \ \ \forall d\in C(x^*),
    \end{equation}
 and by Theorem \ref{dual-conditions}
\begin{eqnarray}
  \Upsilon(x^*,\lambda)(d) &:=& \nabla^2_{xx}L(x^*,\lambda)(d,d)- \sigma\big((\lambda^G,\lambda^H)|T^2_\Omega(G(x^*),H(x^*);\nabla G(x^*)d,\nabla H(x^*)d)\big) \nonumber\\
  &=& d_1^2\geq 0, \ \ \ \forall d\in C(x^*). \label{re-eq-1}
    \end{eqnarray}
   (\ref{re-eq-2}) and (\ref{re-eq-1}) indicate that $\nabla^2_{xx}L(x^*,\lambda)$ is positive semidefinite over $C(x^*)$ while
     $\Upsilon(x^*,\lambda)$ is positive definite over $C(x^*)\backslash\{0\}$. In this example, the second-order necessary conditions  involving the second-order tangent set (\ref{re-eq-1})  is stronger than the one not involving the second-order tangent set (\ref{re-eq-2}).

\section*{Acknowledgments.} {The authors are gratefully indebted to the anonymous referees for their valuable suggestions that helped us greatly improve the original presentation of the paper.}



\begin{thebibliography}{1}
\bibitem{AG} {\sc F. Alizadeh and D. Goldfarb}, {\em Second-order cone programming}, Math. Program.,  95(2003), pp. 3-51.


\bibitem{BCS99} {\sc J.F. Bonnans, R. Cominetti, and  A. Shapiro,} {\em Second order optimality conditions based on parabolic second order tangent sets}, SIAM J.  Optim., 9(1999), pp. 466-492.

 \bibitem{BR05}
 {\sc J.F.\ Bonnans and H.\ Ram\'{i}rez C.},
 {\em Perturbation analysis of second-order cone programming problems},
 Math. Program.,  104(2005), pp. 205-227.


 \bibitem{BS}
 {\sc J.F. Bonnans and A. Shapiro},
 {\em Perturbation Analysis of Optimization Problems}, Springer, 2000.
 \bibitem{CT} {\sc J.-S. Chen and P. Tseng}, {\em An unconstrained smooth minimization reformulation of the second-order cone complementarity problem},  Math.  Program.,  104(2005), pp. 293-327.



\bibitem{Clarke1983}
{\sc F.H. Clarke},
{\em Optimization and Nonsmooth Analysis}, Wiley-Interscience,
New York, 1983.

 \bibitem{Cons06} {\sc E. Constantin}, {\em Second-order necessary conditions based on second-order tangent cones},
     Math. Sci. Res. J.,
   10(2006), pp. 42-56.



 \bibitem{Con16} {\sc E. Constantin}, {\em Second-order necessary conditions for set constrained nonsmooth optimization problems via second-order projective tangent cones}, Lib. Math., 36(2016), pp. 1-24.





%
%
%
%

 \bibitem{FLT02}
 {\sc M.\ Fukushima, Z.-Q.\ Luo, and P.\ Tseng},
 {\em Smoothing functions for second-order cone complementarity problems},
 SIAM J. Optim., 12(2002), pp. 436--460.


 \bibitem{GM11} {\sc H.\ Gfrerer}, {\em First-order and second-order characterizations of metric subregularity and calmness of constraint set mappings}, SIAM J. Optim., 21(2011), pp. 1439-1474.

 \bibitem{GM16} {\sc H.\ Gfrerer and B.S.\ Mordukhovich}, {\em Robinson regularity of parametric constraint systems via variational analysis}, SIAM J.  Optim.,  27(2017), pp. 438-465.

    \bibitem{GO16}
 {\sc H.\ Gfrerer and J.V.\ Outrata},
 {\em On computation of generalized derivatives of the normal-cone mapping and their applications}, Math. Oper. Res.,
 41(2016), pp. 1535-1556.

\bibitem{HYF} {\sc S. Hayashi, N. Yamashita, and M. Fukushima}, {\em Robust Nash equilibria and  second-order cone complementarity problems}, J. Nonlinear  Convex Anal., 6(2005), pp. 283-296.


 \bibitem{GJN10}
 {\sc G.\ Giorgi, B.\ Jimenez, and V.\ Novo},
 {\em An overview of second-order tangent sets and their applications to vector optimization},
 SeMA J., 52(2010), pp. 73-96.

%
%
%
 \bibitem{JLZ15}
 {\sc Y.\ Jiang, Y.J.\ Liu, and L.W.\ Zhang},
 {\em Variational geometry of the complementarity set for second-order cone},
 Set-Valued Var. Anal., 23(2015), pp. 399--414.

 \bibitem{JN04}
 {\sc B.\ Jimenez and V.\ Novo},
 {\em Optimizality conditions in differentiable vector optimization via second-order tangent sets},
 Appl. Math.  Optim., 49(2004), pp. 123--144.



%
%
%
%

 \bibitem{OS08}
 {\sc J.V.\ Outrata and D.F.\ Sun},
 {\em On the coderivative of the projection operator onto the second-order cone},
 Set-Valued Anal., 16(2008), pp. 999--1014.

 \bibitem{RW98}
 {\sc R.T.\ Rockafellar and R.J.\ Wets},
 {\em Variational Analysis}, Springer, New York, 1998.

%



 \bibitem{YZ15}
 {\sc J.J. Ye and J.C. Zhou},
 {\em First-order optimality conditions for mathematical programs with
 second-order cone complementarity constraints},  SIAM J. Optim. 26(2016), pp. 2820-2846.

  \bibitem{yzhou} {\sc J.J. Ye and J.C. Zhou}, {\em Exact formula for the proximal/regular/limiting normal cone of the second-order cone complementarity set},  Math. Program., 162(2017), pp. 33-50.


 \bibitem{YZ16}
 {\sc J.J. Ye and J.C. Zhou}, {\em Verifiable sufficient conditions for the error bound property of second-order cone complementarity problems}, Math. Program., 171(2018), pp. 361-395.

 \bibitem{ZZX13}
 {\sc L.W.\ Zhang, N.\ Zhang, and X.T.\ Xiao},
 {\em On the second-order directional derivatives of singular values of matrices
 and symmetric matrix-valued functions},
 Set-Valued Var. Anal., 21(2013), pp. 557--586.


%

 \bibitem{ZTC15}
 {\sc J.C.\ Zhou, J.Y.\ Tang, and J.-S.\ Chen},
 {\em Parabolic second-order directional differentiability in the Hadamard sense
 of the vector-valued functions associated with circular cones},
 J. Optim. Theory Appl. 172(2017), pp. 802-823.



 \end{thebibliography}
 \end{document}